\documentclass[12pt]{article}
\usepackage{amsmath}
\usepackage{amssymb}
\usepackage{amsthm}
\usepackage{cite}
\usepackage{graphicx,color}

\newcommand{\sH}{{\mathcal H}}

\def\a{{\alpha}}
\def\b{\beta}

\def\de{\Delta}

\def\g{\gamma}
\def\ga{\Gamma}

\def\la{\Lambda}
\def\s{\sigma}
\def\si{\Sigma}

\def\va{\varphi}
\def\va{\varphi}
\def\ve{\varepsilon}

\def\om{\Omega}
\def\th{\theta}

\def\vp{\varepsilon}
\def\vt{\vartheta}

\def\ts{\times}

\def\iy{\infty}
\def\im{{\rm Im\, }}

\def\Up{\Upsilon}

\def\wt{\widetilde}

\def\ov{\overline}
\def\spa{{\rm Span}}

\def\BC{{\mathbb C}}

\def\BN{{\mathbb N}}

\def\BR{{\mathbb R}}

\renewcommand{\theequation}{\arabic{section}.\arabic{equation}}

\newcommand{\bpr}{{\noindent\textbf{Proof.}\ \ }}
\newcommand{\epr}{{$\mbox{}$ \hfill $\Box$}}

\newtheorem{Pa}{Paper}[section]
\newtheorem{Tm}[Pa]{{\bf Theorem}}

\newtheorem{Cy}[Pa]{{\bf Corollary}}
\newtheorem{Rk}[Pa]{{\bf Remark}}

\newtheorem{Pn}[Pa]{{\bf Proposition}}
\newtheorem{prop}[Pa]{Proposition}
\newtheorem{lem}[Pa]{Lemma}
\newtheorem{cor}[Pa]{Corollary}

\newtheorem{defi}[Pa]{Definition}

\newtheorem{remark}[Pa]{Remark}

 \newtheorem{thm}[Pa]{\bf Theorem}

\def\Lam{\Lambda}
\def\va{\varphi}
\def\vp{\varphi}
\def\vt{\vartheta}
\def\mS{\mathcal{S}_{m_2\times m_1}} 

\newcommand{\ands}{\quad\mbox{and}\quad}

\newcommand{\nn}{\nonumber}

\title{{Skew-selfadjoint Dirac systems with  rational}  {rectangular Weyl }
{functions: explicit solutions of direct and inverse problems   and}  {integrable wave equations}}

\author{B. Fritzsche, M.A. Kaashoek, B. Kirstein,  A.L. Sakhnovich}

\date{}

\begin{document}

\maketitle

\begin{abstract} 
In this paper we study {direct and inverse problems for} discrete and continuous time  skew-selfadjoint 
Dirac  systems  with rectangular (possibly non-square) pseudo-exponential potentials. For such a system 
the Weyl function is a strictly proper rational {rectangular} matrix function  and any strictly proper 
rational matrix function appears in this way. In fact, extending earlier results,  given  
a strictly proper rational matrix function we present an explicit procedure  to recover 
the corresponding potential using techniques from mathematical system and control 
theory.  {We also introduce and study a  nonlinear generalized discrete Heisenberg magnet model}, 
extending  earlier  results for the isotropic case. 
A large  part of the paper is devoted to the {related} discrete  time systems of  which the  
pseudo-exponential potential depends on an additional continuous time parameter. 
{Our techniques} allows us to obtain explicit solutions  for the
generalized discrete Heisenberg magnet model {and evolution of the Weyl functions}.

\end{abstract}

{MSC(2010):} 


\section{Introduction} \label{intro}
\setcounter{equation}{0} 

{A skew-selfadjoint Dirac system (also called a pseudo-canonical, Zakharov-Shabat  or AKNS system)} has the form:
\begin{equation}       \label{se1}
\frac{d}{dx}y(x, z )=({i} z j+jV(x))y(x,
z ) \qquad
(x \in \BR_+, \quad z \in \BC),
\end{equation}
\begin{equation}   \label{se2}
j = \left[
\begin{array}{cc}
I_{m_1} & 0 \\ 0 & -I_{m_2}
\end{array}
\right], \hspace{1em} V= \left[\begin{array}{cc}
0&v\\v^{*}&0\end{array}\right],
 \end{equation}
where  $\BC$ stands for the complex plane, $\BR_+$ denotes the non-negative real semi-axis, $I_{m_k}$ is the $m_k \times m_k$ identity matrix, $v(x)$ is an $m_1 \times m_2$ matrix function, which is called the {\textit{potential}} of the system, and $j$ and $V(x)$ are $m\times m$ matrices, $m:=m_1+m_2$.  {Note that $(jV)^*=-jV$, and therefore the system \eqref{se1} is called \textit{skew-selfadjoint}.}

Like {a} selfadjoint Dirac system $\frac{d}{dx}y(x, z )={i} (z j+jV(x))y(x,
z )$, {the} system \eqref{se1} is also an auxiliary system for various important integrable nonlinear wave equations and the case $m_1 \not= m_2$ corresponds to rectangular and multicomponent  versions of these equations. Here we solve explicitly (in terms of Weyl functions) direct and inverse problems for {the system \eqref{se1} for the case when the rectangular (possibly non-square) potential $v$ is  \textit{pseudo-exponential}} (see formula \eqref{7.1}).
 {The direct} problem consists {in constructing} the Weyl function and inverse problem is the problem to recover $v$ from the Weyl function.
 
 We also derive explicit solutions of direct and inverse problems for a discrete analogue of the skew-selfadjoint
 system \eqref{se1}, namely, for the system 
 \begin{align} 
&y_{k+1}(z)=\left(I_m+ i z^{-1} C_k\right) y_k(z),  \quad C_k=U_k^*jU_k,  \  \mbox{where} \nn\\[.15cm]
&\hspace{4.5cm}\ U_k^*U_k=U_kU_k^*=I_m, \quad k=0, 1,2, \ldots. \label{d'1}
\end{align}
This system was studied in \cite{KaS} for the important subcase $m_1=m_2$
and explicit solutions of direct and inverse problems were obtained for that subcase.
See \cite{FKRS14} for  the discrete analogue of selfadjoint Dirac system (and its general Weyl theory).
 
 A large part of the paper is devoted to the case when the system  \eqref{d'1} depends on an additional continuous time parameter.  A special choice of  the additional parameter  allows us to introduce a generalized discrete Heisenberg magnet model. We use our results on system \eqref{d'1} with a general-type $j$ in order to construct explicit  solutions and the evolution of the Weyl function for this generalized model. The results obtained generalize those  for the case $m_1=m_2=1$  in \cite{KaS}  which dealt with the  well-known discrete isotropic Heisenberg magnet model \cite{FT, Skl}.

\medskip
\noindent\textbf{History.} 
 Direct and inverse spectral problems for   Dirac systems have a long and interesting history. For the selfadjoint case  with  a scalar continuous potential $v$  it all started with the seminal paper \cite{Kre1} by M.G. Krein  (see also the discussions in \cite{AGKLS2, SaL50}). This theory is closely related to  the Weyl theory (see, e.g., \cite{AGKLS2, CG2, LS, SaA3, SaSaR, SaL3} and references therein). The Weyl theory for the skew-selfadjoint systems of type   \eqref{se1} was developed much later  \cite{CG2,  FKRS4, GKS2, SaA1, SaA8, SaSaR} and the method of operator identities \cite{SaL1, SaL2, SaL3} played a fundamental role in these studies.

Particular subclasses of general-type potentials $v$ are of special interest, since, for some subclasses, direct and inverse problems may be solved explicitly. Explicit solutions of spectral problems for selfadjoint Dirac systems with the {\it strictly pseudo-exponential\,} square potentials $v$ were given in \cite{AG1, AlGo}. These constructions were based on the procedure to solve  direct and inverse problems for systems with general-type potentials. A few years later direct and inverse problems in terms of spectral and Weyl functions were solved explicitly {in \cite{GKS1}} for a wider class of selfadjoint Dirac systems, that is, for systems with   \textit{pseudo-exponential} square potentials $v$  (see also \cite{GKS6}). The same problem for system \eqref{se1}  was dealt with in \cite{GKS2}. {Moreover,  the above mentioned}  problems for systems with \textit{pseudo-exponential}  potentials were studied directly instead of using procedures for general-type  selfadjoint (or, correspondingly, skew-selfadjoint) Dirac  systems as in \cite{AG1, AlGo}. We note that direct methods in explicit solution of spectral problems go back  to the works by B\"acklund \cite{Ba} and Darboux \cite{Da1, Da2} with essential further developments in \cite{cr, de, fg, gt, kr57}. Finally, direct and inverse problems for selfadjoint Dirac systems with the {\it pseudo-exponential} rectangular {(possibly non-square)} potentials $v$ were solved explicitly  in terms of Weyl functions in \cite{FKRSp1}.

\medskip
\noindent\textbf{Contents.} 
The paper consists of  four sections (the present introduction included) and an appendix. Section \ref{dies} deals with direct and inverse problem for continuous time systems.  In this section we develop  further the results from \cite{GKS2} on  skew-selfadjoint system \eqref{se2} with square potentials, also using  some ideas  and results from the works 
\cite{FKRS4, FKRSp1} on systems with general-type rectangular  (non-square) potentials. Section \ref{ddies}  treats  direct and inverse problems for discrete  time systems, generalizing  earlier results from \cite{KaS} to the rectangular (non-square) case. Finally, Section \ref{NL}  is devoted to generalized discrete Heisenberg magnet model,
which is equivalent to the  compatibility condition of two  systems depending on two parameters $t$ and $k$ where $t$ is non-negative real and $k\in \BN_0$.    For the convenience of the reader  the appendix Section \ref{App} presents   a number of general facts regarding admissible quadruples that are used throughout the paper.  
This theory is closely related to the triples approach \cite{GKS1, GKS2, GKS6, KaS} and generalized B\"acklund-Darboux transformation (GBDT) method \cite{FKRS, FKRSp1, ALS94, SaAJMAA, SaSaR}
(see also references therein).

\medskip
\noindent\textbf{Notation.} We conclude with  some information on notations that are used throughout the paper. As usual $\BN_0$ stands for the set of non-negative integers, i.e., the natural numbers with zero included, and $\BR_+$ denotes the set of non-negative real values. The symbols  $\BC_+\,$  and $\BC_-$  stand  for the upper and lower  half-plane,  respectively. Furthermore, $\ov {\BC_+}$ denotes   the closed upper  half-plane (i.e.,   $\BC_+ \cup \BR$), and $\BC_M$ stands for the  open half-plane $\{z:\, \Im(z)>M>0\}$. By $\|\cdot \|$ we denote the $\ell^2$ vector norm or the induced  matrix norm, $\spa$ stands for the linear span, and $\s(\a)$ stands for the spectrum of $\a$.  The class of $m_2 \times m_1$ contractive matrix functions (Schur matrix functions) 
on some domain $\om$ is denoted by $\mathcal{S}_{m_2 \times m_1}(\om)$.  We write $S>0$ when
the matrix $S$ is positive definite. {The matrix $(\alpha^{-1})^*$ is denoted by $\alpha^{-*}$ and $\im \a$ stands for the image of $\a$.}

\section{Continuous case: direct and inverse \\ problem}\label{dies}
\setcounter{equation}{0}
We begin with introducing the notion of a pseudo-exponential potential. {The starting point is a matrix function $v$ of  the form
\begin{align} \label{7.1}&
v(x)=2\vt _{1}^{*}e^{ix \alpha^{*}} 
S(x)^{-1}e^{i x \alpha }\vt_2, \quad  x\in \BR_+.
\end{align}
Here $\vt_1$ and $\vt_2$ are matrices of sizes  $n\ts m_1$ and $n\ts m_2$, respectively, $\a$ is an $n\ts n$ matrix,   and $S$  is the $n \times n$  matrix function 
given by
\begin{align}   \label{7.2}&
S(x)=S_0+ \int_{0}^{x} \Lambda(t)j \Lambda (t)^{*}dt,
\quad S_0 > 0,
\quad
\la(x):= \begin{bmatrix} e^{-i x \alpha } \vt_{1} 
& e^{i x
\alpha } \vt_2 \end{bmatrix}.
\end{align}
Furthermore, we require the initial value $S_0$  in \eqref{7.2}  to satisfy  the following matrix identity:
\begin{equation} \label{7.3} 
\a S_0-S_0 \a^*=i (\vt_1 \vt_1^*+ \vt_2 \vt_2^*).
\end{equation}}
From \eqref{7.1} and \eqref{7.2} it is clear that  the potential $v$ in \eqref{7.1}  is uniquely determined   by   the  quadruple $\{ \a, S_0,  \vt_1, \vt_2\}$.  In combination with  $S_0$  being  positive definite, the  identity  \eqref{7.3} implies that   $S(x)$  is also positive definite; see \eqref{posSx} below. In particular, $S(x)$ is invertible for each $x\geq 0$, and {hence $v$ is well-defined.}   

{When $S_0>0$ and  \eqref{7.3} holds we call  the quadruple $\{ \a, S_0,  \vt_1, \vt_2\}$   an   \emph{admissible quadruple}.
We call  $v$ in \eqref{7.1} the  \emph{pseudo-exponential potential}
generated by the  admissible quadruple  $\{ \a, S_0,  \vt_1, \vt_2\}$. See appendix Section~\ref{App} for a brief review of properties  of an admissible quadruple and the relation with the theory of S-nodes.}

{The definition of  a pseudo-exponential potential   given here is somewhat different from the definition  in  \cite{GKS2} which starts with \cite[Eq. (0.2)]{GKS2}. However, Proposition 1.1 in  \cite{GKS2}  tells us that  a  pseudo-exponential potential   in  sense of \cite{GKS2} is also a  pseudo-exponential potential as defined above.  With some minor modifications  the reverse implication is also true  (see Proposition \ref{propdefv} at the end of this section).}

In the present  paper, as opposed to  \cite{GKS2}, we do not require the matrices $\vt_1$ and $\vt_2$ to be square, i.e., $m_1$ and  $m_2$ are not required to be equal. Note that in \cite{GKS2} the matrix $S_0$  is just the $n\ts n$ identity matrix  but, as   the next lemma shows, it is convenient to allow $S_0$ to be just positive definite.

\begin{lem}\label{lemadm1}
{Let $\{ \a, S_0,  \vt_1, \vt_2\}$ be an admissible quadruple, and} for each $ x\in \BR_+$  let $\si(x)$ be the   quadruple defined by 
 \begin{equation}
 \label{defsix}
\si(x)=\{ \a, S(x),  e^{-ix\a}\vt_1, e^{ix\a}\vt_2\},
\end{equation}
where $S(x)$ is the matrix defined by the first identity in \eqref{7.2}. {Then  $\si(x)$ is  an  admissible quadruple for  each $ x\in \BR_+$.}
\end{lem}
\bpr   We first  show that
\begin{equation}
\label{7.3a}
\a S(x)-S(x)\a^*=i \la(x)\la(x)^*, \quad  x\in \BR_+.
\end{equation}
To prove this identity note that the relations in \eqref{7.2} are equivalent to
\begin{align}
\frac{d}{dx}\la(x)&=-i\a \la(x)j, \quad \la(0)=\begin{bmatrix}     \vt_{1} 
& \vt_2 \end{bmatrix} \quad ( x\in \BR_+);\label{defLax}\\[.2cm]
\frac{d}{dx} S(x) &=\la(x)j\la(x)^*, \quad S(0)=S_0  \quad (x\in \BR_+). \label{defSx}
\end{align}
Here $j$ is the signature matrix defined \eqref{se2}. If follows that
\begin{align*}
 & \frac{d}{dx} \Big(\a S(x)-S(x)\a^*\Big)=  \a \la(x)j\la(x)^*-  \la(x)j\la(x)^*{\a^*},  \\
   &  \frac{d}{dx} i \la(x)\la(x)^* =\a \la(x)j \la(x)^*- \la(x)j\la(x)^*\a^*.
\end{align*}
On the other hand, using  \eqref{7.3}, the functions $\a S(x)-S(x)\a^*$ and  $ i\la(x)\la(x)^*$  have the same value at $x=0$.  But then \eqref{7.3a} holds true.

It remains to show  that $S(x)$ is positive definite for each $x\geq 0$. To do this we prove the following inequality:
\begin{equation}
\label{posSx}
 e^{-ix\a}S(x) e^{ix\a^*}\geq S_0, \quad x \in \BR_+.
\end{equation}
Using \eqref{defSx} and \eqref{7.3a} we have 
\begin{align}\nn
\frac{d}{dx}\Big(  e^{-ix\a}S(x) e^{ix\a^*}\Big)  & =   e^{-ix\a}\left(\frac{d}{dx}S(x)\right) e^{ix\a^*}\\
&\nn \hspace{2.7cm} -i  e^{-ix\a}\left(\a S(x)-S(x)\a^*\right)e^{ix\a^*}\\[.1cm]
\nn    & = \la(x)j\la(x)^*+ \la(x)\la(x)^*\\[.1cm]
    &= 2e^{-ix\a}\vt_1^*\vt_1e^{ix\a^*}\geq 0, \quad  x\in \BR_+. \nn
\end{align}
But then 
\[
 e^{-ix\a}S(x) e^{ix\a^*}- S_0=\int_0^x \frac{d}{dt}\Big(  e^{-it\a}S(t) e^{it\a^*}\Big) \geq 0, \quad  x\in \BR_+.
\]
In particular,  \eqref{posSx} holds true. \epr

 We note that the quadruple $\{ \a, S_0,  \vt_1, \vt_2\}$ generating $v$ coincides  with $\Sigma(0)$ in \eqref{defsix}.

Using straightforward  modifications of the proof of \cite[Theorem 1.2]{GKS2} (or particular cases of the more  general \cite[Theorem 3]{ALS94} or \cite[Theorem 1.2 and Proposition 1.4]{SaAJMAA}) we obtain the next proposition.

\begin{prop}\label{FundSol}
Let $v$ be  the pseudo-exponential potential  generated by the admissible quadruple  $\si(0)=\{ \a, S_0,  \vt_1, \vt_2\}$, and let $\si(x)$ be the admissible quadruple  defined  by \eqref{defsix}. Then the fundamental solution $u$ of the system \eqref{se1},  normalized by $u(0,z) \equiv I_m$,  where $m=m_1+m_2$, admits the following representation
\begin{align} 
& u(x, z)=W_{\si(x)}(z) e^{i x z j}W_{\si(0)}( z)^{-1},\label{7.4}\\[.2cm]
&W_{\si(x)}(z):=I_m +i \la(x)^*S(x)^{-1}(z I_n- \a)^{-1}\la(x). \label{7.5}
\end{align}
\end{prop}
The \emph{transfer function} of the form \eqref{7.5} (the transfer function in Lev Sakhnovich form) was introduced and studied in \cite{SaL1} (see also \cite{SaL2, SaL3, SaSaR}
and references therein). In our case, we refer to $W_{\si(x)}$  as the the transfer function  associated with the admissible quadruple $\si(x)$; 
see the first paragraph of Section \ref{App}.  The next result is a simple generalization  of \cite[Proposition 1.4]{GKS2}.

\begin{prop}\label{PnBound} 
Pseudo-exponential potentials  $v$ are  bounded on the semiaxis $x \geq 0$.
\end{prop}
\bpr  Assume that the potential $v$ is generated by the admissible quadruple $\si(0)=\{ \a, S_0,  \vt_1, \vt_2\}$, and let $\si(x)$ be the admissible quadruple  defined  by \eqref{defsix}. 
From Lemma \ref{lem11} we know that the eigenvalues of $\a$ belong to $\overline{\BC_+}$. Furthermore, using  the identity  \eqref{basicineq2} with $\si(x)$ in place of $\si$, we see that 
\begin{equation}\label{c2}
i(z-\bar{z})\la(x)^*(\bar{z}I_n-\a^*)^{-1}S(x)^{-1}(zI_n-\a )^{-1}\la(x)\leq I_m,\quad  z\in \BC_-.  
\end{equation}
Since the resolvent $(z I_n-\a)^{-1}$ is well-defined in any open domain
$\om$ in $\BC_-$, we have $\spa_{z\in \om}(z I_n-\a )^{-1}\vt_k \supseteq \vt_k$ for $k=1,2$. Therefore, inequality \eqref{c2} (where $\la(x)$ is given by \eqref{7.2}) implies that
\begin{align} \label{c4}&
\sup_{x \geq 0}\Big(\left\|S(x)^{-1/2}e^{-{i} x
\alpha }\vt_1\right\|+\left\|S(x)^{-1/2}e^{{i} x
\alpha }\vt_2\right\|\Big)<\infty.
\end{align}
From \eqref{7.1} and \eqref{c4}, it  is immediate  that $v$ is bounded, that is, for some $M>0$ we have
 \begin{align}\label{c5}
 \|v(x)\| \leq M  \quad x\in \BR_+,
 \end{align}
which completes the proof.\epr

\paragraph{Weyl function: the direct problem.}  
The concept of  a Weyl function of Dirac system has a long history (see the Introduction). Following {the definition of a Weyl function  for Dirac systems with square potentials (see} {also} \cite{SaSaR} for {the case of non-square potentials}) we say that  {a meromorphic} {function $\vp$ satisfying \eqref{c5} is a \textit{Weyl function} of  the system \eqref{se1}  whenever it satisfies} the inequality
\begin{align}&      \label{c6}
\int_0^{\infty}
\begin{bmatrix}
I_{m_1} & \vp(z)^*
\end{bmatrix}
u(x,z)^*u(x,z)
\begin{bmatrix}
I_{m_1} \\[.2cm] \vp(z)
\end{bmatrix}dx< \infty , \quad z\in \BC_M.
\end{align}
{Here $u$ is the fundamental solution of \eqref{se1} normalized by $u(0,z)=I_m$.} We note that for skew-selfadjoint Dirac systems \eqref{se1}   Weyl functions have been introduced in \cite{FKRS4} in an equivalent but different way.  However, Proposition 2.2 and Corollary 2.8 from \cite{FKRS4} immediately yield the following result concerning direct problem.

\begin{prop}\label{PnC2} 
Let \eqref{c5} hold. Then there is a unique function $\vp$ such that \eqref{c6} is valid.
This function $\vp$ is analytic and contractive in $\BC_M$ $($i.e., $\vp\in \mS(\BC_M))$.
\end{prop}

For the case of   {a pseudo-exponential potential} $v$ we produce an explicit expression for the Weyl function.

\begin{thm}\label{PnC3} Consider the  skew-selfadjoint Dirac system \eqref{se1}, and assume that $v$ is  a  pseudo-exponential potential  generated  by the admissible quadruple $\si(0)=\{\a, S_0, \vt_1,  \vt_2\}$. Then the Weyl function $\vp$ of  the Dirac system \eqref{se1} is given by
\begin{equation} \label{7.9} 
\vp(z)={i} \vt_2^*S_0^{-1}(z I_n
- \a^\ts)^{-1}
\vt_1 , \quad \a^\ts:=\a-{i} \vt_1\vt_1^*S_0^{-1}.
\end{equation}
\end{thm}

\bpr  From the proof of Lemma \ref{lem12} we know that    \eqref{7.9} can be rewritten in the form:
\begin{equation}  \label{c7} 
\vp(z)={i} \vt_2^*S_0^{-1}(zI_n-\a)^{-1}\vt_1\big(I_{m_1}+{i}\vt_1^*S_0^{-1}(z I_n-\a)^{-1}\vt_1\big)^{-1}.
\end{equation}
Taking into account the equivalence of  \eqref{7.9} and  \eqref{c7} and using  \eqref{7.5} together with the second equality in \eqref{7.2} (both at $x=0$), we derive that $\vp$ of the form \eqref{7.9}
satisfies the relation:
\begin{equation} \label{c11} 
\begin{bmatrix}
I_{m_1} \\ \vp(z)
\end{bmatrix}=W_{\si(0)}(z) \begin{bmatrix}
I_{m_1} \\ 0
\end{bmatrix}\big(I_{m_1}+{i}\vt_1^*S_0^{-1}(z I_n-\a)^{-1}\vt_1\big)^{-1}.
\end{equation}
Formulas \eqref{7.4} and \eqref{c11} imply that
\begin{align} \label{c12} 
u(x,z)\begin{bmatrix}
I_{m_1} \\ \vp(z)
\end{bmatrix}=e^{i  x z}W_{\si(x)}(z)  \begin{bmatrix}
I_{m_1} \\ 0
\end{bmatrix}\big(I_{m_1}+i\vt_1^*S_0^{-1}(z I_n-\a)^{-1}\vt_1\big)^{-1}.
\end{align}
Using   $S_0 > 0$  and taking inverses in \eqref{posSx}, we see that 
\[
\sup_{x\geq 0}\|e^{-i x \a^*}S(x)^{-1}e^{i x \a}\|<\infty.
\] 
This yields  the inequality
\begin{align} \label{c13} 
\sup_{x\geq 0}\|e^{-\eta  x }S(x)^{-1}\|<\infty \  \ \mbox{for sufficiently large values of $\eta$}.
\end{align}
From  \eqref{c13} and the identity  \eqref{basicineq2} with $\si(x)$ in place of $\si$, it then follows that the function  { $x\to e^{i (z-\ov z)  x }W_{\si(x)}( z)^*W_{\si(x)}( z) $
belongs to $ L^2_{m\times m}(0, \, \infty)$  for sufficiently large values of $\Im(z)$:
\begin{equation}\label{c14} 
 e^{i (z-\ov z)  x }W_{\si(x)}( z)^*W_{\si(x)}( z) \in  L^2_{m\times m}(0, \, \infty).
\end{equation}
That is,} the entries of $e^{i (z-\ov z)  x }W_{\si(x)}( z)^*W_{\si(x)}( z)$ are squarely summable  with respect to $x$  for sufficiently large values of $\Im(z)$.  Finally, in view of  \eqref{c12} and \eqref{c14}, the inequality \eqref{c6} holds for $\vp$ given by \eqref{7.9} and for sufficiently large values of $\Im(z)$. 

Taking into account Proposition \ref{PnC2} and the analiticity of $\vp$ given by \eqref{7.9}, we see that \eqref{c6} holds for all $z\in \big(\BC_M \setminus \s(\a^\ts) \big)$, that is, $\vp$ is the Weyl function.
\epr

\paragraph{Weyl function: the inverse problem.} Theorem \ref{PnC3} presents the solution for  the direct problem. We now turn to the inverse problem. The uniqueness theorem below is immediate from \cite[Theorem 3.21 and Corollary 3.25]{SaSaR}
(see also \cite{FKRS4}).
\begin{thm}\label{TmC3} Let  $\vp$ be a Weyl function of  a skew-selfadjoint Dirac system with  a   potential $v$ which is bounded on $[0, \, \infty)$.
Then this $v$ can  {be} uniquely recovered from $\vp$.
\end{thm}

For the case of  pseudo-exponential potentials  we have an explicit procedure to recover the potential from the Weyl function.  
 {This procedure uses such well-known notions from {\it control theory} as {\it realization, minimal realization} and {\it McMillan degree} (see, e.g., \cite{BGK1, KFA} or  \cite[App. B]{SaSaR}).
We note}
that, according to Theorem \ref{PnC3},  the Weyl function $\va$  of a  skew-selfadjoint Dirac system with a pseudo-exponential potential is a strictly proper rational matrix function. 

\begin{thm} \label{TmC5}
Let $ \vp $ be a strictly proper rational $m_2 \times m_1$ matrix function. Then $\vp$ is the Weyl function of a skew-selfadjoint Dirac system with  a pseudo-exponential potential $v$
{ generated  by an admissible quadruple.
The corresponding quadruple $ \{\a, S_0,  \vt_1,   \vt_2\}$ can be obtained explicitly by using the  following procedure.
First, we construct a minimal realization  of $ \vp:$
\begin{equation}\label{c15!} 
\va(z)=i \th_2^*(z I_n
- \g)^{-1}
\th_1.
\end{equation}
Next, we choose $X$ to be the unique positive definite solution $X$  of the Riccati equation
$\g X-X\g^*-iX\th_2\th_2^*X +i \th_1\th_1^*=0$.
Finally, we put
\begin{equation} \label{c17a!}
S_0=I_n, \quad \vt_{1}=X^{ -{1}/{2}} \theta_{1}, \quad \vt_{2}=X^{{1}/{2}}  \theta_{2}, \quad
 \a=X^{- {1}/{2}}\g  X^{ {1}/{2}}+i \vt_1\vt_1^*. 
\end{equation} }
\end{thm}
\bpr  {Given a  strictly proper rational  matrix function $\va$, Theorem \ref{thminv} provides} {the described above
three step procedure to construct a quadruple $\{\a, S_0, \vt_1, \vt_2\}$ such that representation \eqref{7.9} of $\vp$ holds.
It also follows from Theorem \ref{thminv} that the procedure is well-defined and the  quadruple $\{\a, S_0, \vt_1, \vt_2\}$ is admissible.
Then we know from Theorem \ref{PnC3} that  $\va$ we started with is precisely  the Weyl function of the skew-selfadjoint Dirac system  with 
the potential $v$  generated by    $\{\a, S_0, \vt_1, \vt_2\}$. }\epr

{\paragraph{A final remark.} } As  {we} note in the second paragraph before Lemma \ref{lemadm1} the definition of a pseudo-exponential  {potential  given} in the beginning of this section differs from the definition employed in \cite{GKS2} which starts from formula (0.2) in \cite{GKS2}. The next  proposition presents the analogue of formula (0.2) in \cite{GKS2},
{which coincides with (0.2) for the case $m_1=m_2$ and $S_0=I_n$.}
The conclusion is that the two definitions lead to the same class of potentials.

{\begin{prop}\label{propdefv} Let $v$ be the pseudo-exponential potential   generated by the admissible quadruple  $\{ \a, S_0,  \vt_1, \vt_2\}$, and let $A$ be the $2n\ts 2n$ matrix defined by
\begin{equation}
\label{defA}
A=\begin{bmatrix} \a^* &0\\ -\vt_1 \vt_1^*&\a \end{bmatrix}.
\end{equation}
Then the potential $v$ is also given by 
\begin{equation}
\label{altdefv}
v(x)=2 \vt_1^* \Big(\begin{bmatrix}S_0 & -iI_n\end{bmatrix}e^{-2ixA}\begin{bmatrix}I_n \\ 0\end{bmatrix}\Big)^{-1}\vt_2, \quad x\in \BR_+.
\end{equation}
\end{prop}}
 {The proof of the proposition is given in Appendix (and  is close to the considerations in \cite{GKS2}).}

\section{Discrete case: direct  and inverse problem} \label{ddies}
\setcounter{equation}{0}

Recall that the discrete skew-selfadjoint Dirac system  {(SkDDS)}   has the form:
\begin{equation} \label{d1}
y_{k+1}(z)=\left(I_m+ i z^{-1}
C_k\right)
y_k(z),  \quad C_k=U_k^*jU_k \quad \left( k \in \BN_0
\right).
\end{equation}
Here $\BN_0$ stands for the set of non-negative integers, the  matrices $U_k$ are unitary of size $m \ts m$, where $m$ does not depent on $k$,  and   $j$ is the $m \ts m$ signature matrix defined by \eqref{se2} with $m=m_1+m_2$ and with $m_1$ and $m_2$ not depending on $k$.  The sequence  $\{C_k\}_{k\in \BN_0} $ is called the \emph{potential} of the system.  Note that the second part of \eqref{d1} implies that $C_k=C_k^*= C_k^{-1}$ for each $k\in \BN_0$.

\begin{defi}\label{DnD0} The Weyl function of SkDDS is an $m_1 \times m_2$ {meromorphic} matrix function $\vp(z)$ on $\BC_M$ $($for some $M>0)$,
which satisfies the inequality
\begin{equation} \label{d1!}
\sum_{k=0}^{\infty}\begin{bmatrix}  \vp(z)^* & I_{m_2}\end{bmatrix}w_k(z)^* w_k(z)
\left[
\begin{array}{c}
 \vp(z) \\ I_{m_2}
\end{array}
\right] < \infty,
\end{equation}
where $w_k(z)$ is the {fundamental solution} of SkDDS normalized by $w_0(z) \equiv I_m$.
\end{defi}

We shall consider the case when   the potentials $\{C_k\}_{k\in \BN_0} $ are \emph{pseudo-exponential} (see Definition \ref{DnPE} below), and we shall show that for  such potentials  the Weyl function can be constructed explicitly.  

Pseudo-exponential potentials have been  introduced in \cite{KaS} for  the square case when $m_1=m_2$.  We shall show  {that the  same} scheme developed in \cite{KaS} {for  constructing} the corresponding Weyl  function  {also  works in} the non-square case when  $m_1$  and  $m_2$ are not equal.

 {Similar to the continuous case},  {our starting point} to define pseudo-expo-nential potentials is an admissible  {quadruple $\si_0=\{ \a, S_0,  \vt_1, \vt_2\}$, and $\la_0= \begin{bmatrix} \vt_1&\vt_2 \end{bmatrix}$}. In addition, we assume that $\a$ is non-singular.   {We set}
\begin{align}
   \la_{k+1}&= \la_k+i\a^{-1}\la_k j   \quad ({k \in \BN_0}); \label{defLak}\\[.2cm]
   S_{k+1} & =  S_{k} +\a^{-1}S_k \a^{-*}+ \a^{-1}\la_kj\la_k^* \a^{-*} \quad ({k \in \BN_0}),  \label{defSk}
\end{align}
 {where $\a^{-*}$ stands for $(\a^{-1})^*$.}
From \eqref{defLak} and   $\la_0= \begin{bmatrix} \vt_1&\vt_2 \end{bmatrix}$ it follows that 
\begin{equation}
\label{d21}
\la_k=\begin{bmatrix}(I_n+i \a^{-1})^{k}\vt_1 & (I_n-i \a^{-1})^{k}\vt_2 \end{bmatrix}  \quad ({k \in \BN_0}).
\end{equation}
By definition $\si_k$ is the quadruple  given by
\begin{equation}
\label{defsik} 
\si_k= \{\a, S_k, (I_n+i \a^{-1})^{k}\vt_1,  (I_n-i \a^{-1})^{k}\vt_2\}.
\end{equation} 
The next lemma shows that this quadruple is again admissible if, in addition, the pair $\{\a, \vt_1\}$ is controllable.  

\begin{lem}\label{lemadmdiscrk}
 Let   {$\si_0=\{ \a, S_0,  \vt_1, \vt_2\}$} be an admissible quadruple.  Assume that   the pair $\{\a, \vt_1\}$ is controllable.   Then $\a$ is non-singular,  and the quadruple  $\si_k$ defined by \eqref{defsik} is admissible for each $k \in \BN_0$.  Moreover the pair $\{\a,  (I_n+i \a^{-1})^{k}\vt_1\}$ is controllable. 
\end{lem}

\bpr  The fact that  the pair $\{\a, \vt_1\}$ is controllable  implies that  the pair  $\{\a, \la_0\}$ is also controllable. But then   item (ii) in Lemma \ref{lem11} tells us that $\s(\a)\subset \BC_+$.  In particular, $\a$ is non-singular, and thus the quadruple  $\si_k$  is well-defined. 

Next we use \eqref{defLak} and \eqref{defSk} to prove that 
\begin{equation} \label{d2'}
\a S_k-S_k \a^*= i\la_k \la_k^*  \quad ({k \in \BN_0}).
\end{equation}
This will be done  by induction. Since $\la_0= \begin{bmatrix} \vt_1&\vt_2 \end{bmatrix}$ and the quadruple $\si_0=\{ \a, S_0,  \vt_1, \vt_2\}$ is admissible, the identity \eqref{d2'} holds for $k=0$.  Suppose now that (\ref{d2'}) holds for $k=r$. Then, using the expression for
$S_{r+1}$ from \eqref{defSk}, we obtain 
\begin{align} 
\a S_{r+1}- S_{r+1} \a^* &= i \la_r \la_r ^* + i\a^{-1} \la_r \la_r ^* \a^{-*}+\nonumber \\[.2cm]
&\hspace{1.5cm}+\la_r  j \la_r ^* \a^{-*}- \a^{-1}\la_r j\la_r ^*. \label{d8a}
\end{align}
On the other hand, using \eqref{defLak} for $k=r$, we obtain
\begin{align} 
i\la_{r+1}\la_{r+1}^*&=i(\la_r+i\a^{-1}\la_r j)(\la_r^*-i j \la_r^*\a^{-*})\nonumber \\[.2cm]
&=i \la_r \la_r ^*- \a^{-1}\la_r j\la_r ^*+\la_r  j \la_r ^* \a^{-*}+ i\a^{-1} \la_r \la_r ^* \a^{-*}.\label{d8b}
\end{align}
Together \eqref{d8a} and  \eqref{d8b} yield \eqref{d2'} for $k=r+1$,  and thus \eqref{d2'} holds for all $k\in \BN_0$.

Next we show that $S_k$ is positive definite for each $k\in \BN_0$. To do this we first use that $\s(\a)\subset \BC_+$. The latter inclusion and the identities   \eqref{d2'} imply (cf.,  {\cite{IvSa}
or} Theorem I.4.1 in \cite{GGK1}; also \eqref{altS0})  that
\begin{equation} \label{d23}
S_k=\frac{1}{2\pi}\int_{-\infty}^{\infty}(\lambda  I_n- \a)^{-1}\la_k\la_k^*(\lambda  I_n-\a^*)^{-1}d\lambda , \quad k\in \BN_0.
\end{equation}
Now, the inequality $S_k > 0$ is proved by contradiction.  Assuming that $S_k\not  >0$, we derive (from \eqref{d23}) the existence of a  vector $g\not=0$ such that
$g(\lambda  I_n- \a)^{-1}\la_k\equiv 0$. In particular, using the identitiy  \eqref{d21}, we see that $\wt{g} (\lambda  I_n- \a)^{-1} \vt_1\equiv 0$, where $\wt{g}=g(I_n+{i} \a^{-1})^{k}$. The fact that $\s(\a)\subset \BC_+$ implies that $-i$ is not an eigenvalue of $\a$, and hence the matrix $I_n+{i} \a^{-1}$ is non-singular. But then $\wt{g}\not = 0$ because $g\not =0$.   {On the other hand}, from the Taylor expansion of $\wt g(\lambda  I_n- \a)^{-1}\vt_1$ we obtain that
\[
\wt g \a^p \vt_1=0 \quad (\wt g \not=0), \quad \mbox{for all $p\geq 0$},
\]
which  contradicts  the assumption that the pair $\{\a, \, \vt_1\}$ is controllable.  Thus $S_k$ is positive definite, and hence $\si_k$ is admissible. 

It remains to prove the pair $\{\a,  (I_n+i \a^{-1})^{k}\vt_1\}$ is controllable. To do this note that
\begin{equation}\label{contr1a} 
\spa \bigcup_{\nu=0}^{n-1}\im \a^\nu (I_n+i \a^{-1})^{k}\vt_1 =(I_n+i \a^{-1})^{k}\Big(\spa \bigcup_{\nu=0}^{n-1}\im \a^\nu\vt_1\Big).
\end{equation}
Since the pair $\{\a, \vt_1\}$ is controllable,  the space $\spa \bigcup_{\nu=0}^{n-1}\im \a^\nu\vt_1$  is equal to $\BC^n$. It follows that the space in the right-hand side of \eqref{contr1a} is also equal to $\BC^n$. But then the same holds true  for the space in the left-hand side of \eqref{contr1a}, which implies that  $\{\a,  (I_n+i \a^{-1})^{k}\vt_1\}$ is controllable.
\epr

\medskip
Note that the conditions in Lemma \ref{lemadmdiscrk}  imply that the matrix $S_k$ is non-singular  for each $k$. This allows us to  define the following sequence of matrices:
\begin{equation} \label{d5}
C_k=j+ \la_k^* S_k^{-1} \la_k - \la_{k+1}^* S_{k+1}^{-1} \la_{k+1},
\quad k=0,1,2,\ldots.
\end{equation}

\begin{defi}\label{DnPE} We call the  quadruple  $\si_0=\{ \a, S_0,  \vt_1, \vt_2\}$ strongly  admissible if it is admissible and the pair $\{a, \vt_1\}$ is controllable.  In this case we  refer to  the sequence  of matrices $\{C_k\}_{k\in \BN_0}$ in \eqref{d5} as  the  pseudo-exponential potential generated by  $\si_0$.
\end{defi}
Later we shall  see (Proposition \ref{PnDa1} below)   that the second part of \eqref{d1} is  fulfilled for any pseudo-exponential potential $\{C_k\}_{k\in \BN_0}$.  

\begin{Rk}\label{Rkv1} {We note that for any continuous pseudo-exponential potential $v$ $($see its definition at the beginning of Section \ref{dies}$)$ there is a strongly  admissible
quadruple which generates it. Indeed, recall that in view of Proposition \ref{PnBound} and Theorem \ref{TmC3} there is a unique solution of the inverse problem
considered in Theorem \ref{TmC5}. Thus, starting from the Weyl function of a  system with the pseudo-exponential $v$ we recover some admissible quadruple
generating this $v$ via formula \eqref{c17a!} from Theorem \ref{TmC5}.
It is immediate from \eqref{c17a!} that the corresponding pair $\{a, \vt_1\}$ is controllable, and so the recovered quadruple is strongly  admissible.}
\end{Rk}

Next, we present   a discrete analog of formula \eqref{7.4}  {for} the fundamental solution. See the next  theorem which is a minor generalization  of Theorem 0.1 in \cite{KaS}.

\begin{thm} \label{TmD1}
 Let  $\si_0=\{ \a, S_0,  \vt_1, \vt_2\}$ be a strongly   admissible quadruple, and let $\{ C_k\}_{k\in \BN_0} $ be the pseudo-exponential {potential} generated  by $\si_0$. Then  the fundamental solution $w_k(z )$ $(k\in \BN_0)$ of the normalized initial value problem 
\begin{equation} \label{d1!!}
w_{k+1}(z)=\left(I_m+ \frac{i}{ z}
C_k\right)
w_k(z), \quad w_0(z)\equiv I_m
\end{equation}
can be represented in the form
\begin{equation}\label{d6}
w_k( z )=W_{\si_{k}}( -z ) \left(
I_{m}+\frac{i}{z}j \right)^{k} W_{\si_{0}}( - z)^{-1},
\end{equation}
where  {$\si_k$ is   the admissible quadruple defined by \eqref{defsik} and $W_{\si_{k}}$ is the transfer  function associated with $\si_k$},  that is, 
\begin{equation}\label{d7}
W_{\si_{k}}(\lambda)=I_{m}+i \la_k^{*}
S_k^{-1} ( \lambda I_{n}- \alpha )^{-1} \la_k.
\end{equation}
\end{thm}
\noindent
See the first paragraph of Section \ref{App} for the definition of  the  {transfer function}  associated with an admissible quadruple. 

\smallskip
\bpr  The proof is based on  the following identity:
\begin{equation} \label{d9}
W_{\si_{k+1}} (\lambda)\left(I_{m}-
\frac{i}{\lambda}j\right)=\left(I_{m}- \frac{i}{\lambda }C_k\right)W_{\si_{k}} (\lambda).
\end{equation}
This identity is the analogue of formula (2.2) in \cite{KaS}, and its  proof  follows the same line of reasoning as  the proof of  formula (2.2) in \cite{KaS}. We omit the details. 

 {Next, using \eqref{d9}, we complete the proof  by induction.}  First notice that the equality \eqref{d6} holds for $k=0$.  Next, assume \eqref{d6} is proved for $k=r$. Then,    using   \eqref{d1!!}, 
 \eqref{d6} with $k=r$, and \eqref{d9} with $k=r$ and $\lambda=-z$, we see that
\begin{align*}
w_{r+1}&= \left(I_m+ \frac{i}{ z} C_r\right) w_r(z)\\ 
&=\left(I_m+ \frac{i}{ z} C_r\right) W_{\si_{r}}(-z)\left(
I_{m}+\frac{i}{z}j \right)^{r} W_{\si_{0}}( - z)^{-1}\\
&= W_{\si_{r+1}}(-z)\left(I_m+\frac{i}{z}j\right) \left(
I_{m}+\frac{i}{z}j \right)^{r} W_{\si_{0}}( - z)^{-1}\\
&=W_{\si_{r+1}}(-z)\left(
I_{m}+\frac{i}{z}j \right)^{r+1} W_{\si_{0}}( - z)^{-1}.
\end{align*}
 {Thus \eqref{d6} holds for $k=r+1$, and so it holds for all $k\in \BN_0$.} \epr

 {Our next proposition shows that  the  second equality   in \eqref{d1} is valid.}

\begin{prop} \label{PnDa1} Let  $\si_0=\{ \a, S_0,  \vt_1, \vt_2\}$ be a strongly admissible quadruple, and let $\{ C_k\}_{k\in \BN_0} $ be the pseudo-exponential  {potential} generated  by $\si_0$.  Then we  have
\begin{align} \label{da0} &
C_k=U_k^*jU_k , \quad {\mathrm{where}} \quad U_k^*U_k=I_m  \quad(k\in \BN_0).
\end{align}
\end{prop}

\bpr   {According to Lemma \ref{lemadmdiscrk}, $\si_k$ is an admissible quadruple, and so we may substitute $\si_k$ (instead of $\si$) into
formula \eqref{basicineq1} from Lemma \ref{lem11}.
 Then  \eqref{basicineq1}  implies that }
\begin{equation} \label{d19}
W_{\si_{k}}(\lambda)W_{\si_{k}}(\overline{\lambda})^* =I_{m} \quad (k \in \BN_0),
\end{equation}
where $\overline{ \lambda}$ stands for complex conjugate for $\lambda$.  Note that \eqref{d19} makes sense only when $\lambda$ and $\overline{\lambda}$  are no poles of  $W_{\si_k}$. 

It is immediate from (\ref{d5}) that $C_k=C_k^*$.  From \eqref{d9} and \eqref{d19} we see that 
\[ 
\left(I_{m}- i\lambda^{-1}C_k\right) \left(I_m+ i\lambda^{-1}C_k\right)=\lambda^{-2}(\lambda^2+1)I_{m}, \quad \lambda\in \BR,
\]
that is, $C_k^2=I_m$. Thus, the equality $C_k=C_k^{-1}$  is valid, and we obtain the representation 
\begin{align} \label{d5!} &
C_k=\wt U_k j_k\wt U_k^*,
\end{align}
where $\wt U_k$ are unitary matrices and $j_k$ are diagonal matrices, the entries of which take the values $\pm 1$.

Next we show that we may choose $j_k=j$.  From Lemma \ref{lem11} we know that that $\s(\a)\subset  \overline{\BC_+}$. In particular,  $-i $ does not belong to  $\s(\a)$. We first deal with the case when   the same holds true for $ i$. In that case \eqref{d19} makes sense for both $\lambda=i$ and $\lambda=-i$.  Let us partition the transfer matrix function  {$W_{\si_{k}}(\lambda)$}  into the two blocks, the first consisting of  the first  $m_1$ columns and the second of the remaining $m_2$ columns:  
\[
W_{\si_{k}}(\lambda)=\begin{bmatrix}  \big(W_{\si_{k}}(\lambda) \big)_1 & \big( W_{\si_{k}}(\lambda))_2 \end{bmatrix}.
\]
Now, take $\lambda =-i$ in \eqref{d9} and multiply the resulting identity from the right by $W_{\si_k}(i)^*$. Then, using \eqref{d19} with $\lambda=-i$, we obtain
\begin{equation} \label{da1}
2\big(W_{\si_{k+1}}( -i) \big)_1 \big(W_{\si_{k}}( i)\big)_1^*=I_m+C_k.
\end{equation}
Repeating the argument with   $\lambda =i$ in place of $\lambda  =-i$ we see that 
\begin{equation} \label{da2}
2\big(W_{\si_{k+1}}(i) \big)_2 \big(W_{\si_{k}}(-i) \big)_2^*=I_m-C_k.
\end{equation}
Formulas \eqref{da1} and \eqref{da2} imply that
\begin{align} \label{da3} &
{\mathrm{rank}}( I_m+C_k)\leq m_1, \quad {\mathrm{rank}}( I_m-C_k)\leq m_2.
\end{align}
The existence of   the representation \eqref{da0} now easily follows from \eqref{d5!} and \eqref{da3}. Thus the proposition is proved when  $ i \not\in \s(\a)$.

It remains to consider the case  when $i$ belongs to $\s(\a)$.  In this case  we   approximate the   original quadruple $\si_0=\{\a, S_0, \vt_1, \vt_2\}$  by  a set of new quadruples $\si(\ve)=\{\a+\ve I_n, S_0, \vt_1, \vt_2\}$, with   $\ve =\ov{\ve}>0$  and $\ve$ sufficiently small. These new quadruples are again strongly admissible.  Thus they   satisfy the  conditions of  Theorem \ref{TmD1} and the additional condition 
\[ i \not\in \s(\a+\ve I_n).
\] 
Therefore, applying the result of the previous paragraph, \eqref{da0} is valid also in this case. Taking limits for $\ve \downarrow 0$ one obtains \eqref{da0} for the original quadruple $\si_0$. \epr

\begin{remark} \label{RkGen} It is not difficult  to see that Theorem \ref{TmD1}, the identities \eqref{d2'} and  the representations in \eqref{da0} hold under weaker  conditions. Indeed, assume that   {$\a$ is non-singular, {that} $S_0=S_0^*$, {that} \eqref{d2'} holds for {$k=0$, and that} $\la_{k}$ and $S_k$ are given by \eqref{defLak} and \eqref{defSk}, respectively. Then the identities \eqref{d2'} are valid for $k\in \BN_0$.} Assuming additionally that $\det \, S_{r} \not=0$  for $0 \leq r \leq N$,  one can show that the fundamental solution $w_k(z )$ $(0 \leq k \leq N)$ of the discrete system  \eqref{d1!!} can be represented in the form
\eqref{d6} and that \eqref{da0} is also valid for  $0 \leq k
\leq N-1$.
\end{remark}

\paragraph{Weyl function: the direct problem.} The following result is a discrete analogue of  Theorem \ref{PnC3}. 
 {The special case $m_1=m_2$ of this result} 
is a somewhat stronger version of Theorem 0.4 in \cite{KaS}.

\begin{thm} \label{TmD5}  Let  $\si_0=\{ \a, S_0,  \vt_1, \vt_2\}$ be a strongly admissible quadruple, and let $\{ C_k\}_{k\in \BN_0} $ be the pseudo-exponential sequence generated  by $\si_0$.  Then the discrete skew-selfadjoint Dirac system
\eqref{d1} has a unique Weyl function $\varphi$ which  is given by the formula
\begin{equation} \label{d24}
\varphi ( z )=-{i} \vartheta_{1}^{*}S_0^{-1}(z  I_{n} + \beta )^{-1}
\vartheta_{2}, \quad \beta:=\a-{i} \vt_2\vt_2^*S_0^{-1}.
\end{equation}
Moreover, $\va$ satisfies \eqref{d1!} in the half-plane $\Im z > \frac{1}{2}$
$($a finite number of points excluded$)$.
\end{thm}
\bpr   \textsc{Part 1.}
 {In order to show that $\va$ is a Weyl function {one first compares} \eqref{d24} with \eqref{defWeyl2} to see that $\va$ is the (second) {function associated with $\si_0$.} 
Thus, {the second part in \eqref{abcd}}
implies that 
\begin{equation} \label{d27}
\vp(z)= b( -z )d( -z )^{-1},
\end{equation}
where $b$ and $d$ are the blocks of $W_{\si_0}$ (see \eqref{d26}).} Using \eqref{d27} we have
\begin{align}\label{d28}&
W_{\si_0}(-z)^{-1}\begin{bmatrix} \vp(z) \\ I_{m_2} \end{bmatrix}
=W_{\si_0}(-z)^{-1}\begin{bmatrix} b(-z) \\ d(-z) \end{bmatrix}d(-z)^{-1}
=\begin{bmatrix} 0 \\ I_{m_2}  \end{bmatrix}d(-z)^{-1}.
\end{align}
Next, for $k\in \BN_0$,  let $W_{\si_k}$ be the transfer function of the quadruple $\si_k$  defined by \eqref{defsik}. Recall that the normalized fundamental solution $w_k$ of \eqref{d1} is given by \eqref{d6}.  Using  \eqref{d28} this yields
\begin{align}
w_k(z)\begin{bmatrix} \vp(z) \\ I_{m_2} \end{bmatrix}&=W_{\si_{k}}( -z ) \left(
I_{m}+\frac{i}{z}j \right)^{k} W_{\si_{0}}( - z)^{-1}\begin{bmatrix} \vp(z) \\ I_{m_2} \end{bmatrix}\nonumber \\
&=W_{\si_{k}}( -z ) \left(
I_{m}+\frac{i}{z}j \right)^{k}\begin{bmatrix} 0 \\ I_{m_2}  \end{bmatrix}d(-z)^{-1}\nonumber \\
&=\left(\frac{z-{i}}{z}\right)^k W_{\si_k}( -z)\begin{bmatrix} 0 \\ I_{m_2}  \end{bmatrix}d(-z)^{-1}.\label{d29}
\end{align}
Since $S_k > 0$ and $\s(\a)\subset \overline{\BC_+}$ (according to item (i) in Lemma \ref{lemadmdiscrk}), the  identity \eqref{basicineq2} implies that
\begin{equation}\label{d25!} 
W_{\si_k}( -z)^{*}
W_{\si_k}( -z )  \leq I_m \quad (z\in \BC_+).
\end{equation}
From \eqref{d29} and \eqref{d25!} we easily derive \eqref{d1!} for all $z\in \BC_{1/2}$,
a finite number of eigenvalues excluded.

\smallskip
\noindent
\textsc{Part 2.} 
It remains to prove the uniqueness of the Weyl function.   {We shall again need the identity \eqref{basicineq2} (for the {quadruple} $\si_k$):}
\begin{equation}\label{d25} 
W_{\si_k}(\lambda)^* W_{\si_k}(\lambda)= 
I_{m}-i( \lambda - \overline{ \lambda})
\la_k^{*}( \overline{ \lambda}I_{n}- \a ^{*})^{-1} S_k^{-1} (  \lambda  I_{n}-  \alpha )^{-1} \la_k.
\end{equation}

 {Since the left upper block of \eqref{d25} is nonnegative,  taking into account
\eqref{d21},
we  obtain for $z=(-\lambda) \in \BC_+$} the inequality
\begin{align}
& \vt_1^{*}( \overline{ z }I_{n}+ \alpha
^{*})^{-1}
(I_n-{i}( \a^*)^{-1})^k S_k^{-1} 
(I_n+{i} \a^{-1})^k(  z I_{n}+  \alpha )^{-1} \vt_1 \leq
\frac {i}{z-  \overline{ z }} I_{m_1}. \label{d33}
\end{align}
Using the same arguments as at the end of the proof of Lemma  \ref{lemadmdiscrk},
we see that  the controllability of the pair $\{\a, \, \vt_1\}$ implies the equality
\begin{align}&
\spa_{z\in O_{\ve}(z_0)}(  z I_{n}+  \alpha )^{-1} \vt_1=\BC^n \label{d34}
\end{align}
for any $\ve$-neighborhood $O_{\ve}(z_0)$ of any $z_0\in \BC$. From \eqref{d33} and \eqref{d34} we see that for  {all} $k\in \BN_0$ the positive definite matrix 
$(I_n-{i}( \a^*)^{-1})^k S_k^{-1}  (I_n+{i} \a^{-1})^k$ is uniformly bounded. Hence,   {the matrix function
\begin{align}\nn
& i(\overline{ z }-z)\vt_1^{*}( \overline{ z }I_{n}+ \alpha
^{*})^{-1}
(I_n-{i}( \a^*)^{-1})^k S_k^{-1} 
(I_n+{i} \a^{-1})^k(  z I_{n}+  \alpha )^{-1} \vt_1 
\end{align}
is sufficiently small for all $z \in \BC_M$ and sufficiently large values of $M$. Recall that the matrix function above coincides 
with the left upper block of the second term on the right hand side of \eqref{d25} (for $\lambda = -z$). Therefore,
relations
\eqref{d6}}  and \eqref{d25} yield the inequality
 {
\begin{align}&
\begin{bmatrix} I_{m_1} & 0\end{bmatrix} W_{\si_0}(-z)^*w_k(z)^*w_k(z)W_{\si_0}(-z)
\begin{bmatrix} I_{m_1} \\ 0\end{bmatrix}\geq \frac{1}{2} I_{m_1}, 
\label{d35}
\end{align}
for $z \in \BC_M$ and sufficiently large values of $M$.} Thus, for sufficiently large $M$ and all $z\in \BC_M$, we have 
\begin{align}&
\dim L \leq m_2, \quad L:=\{g\in \BC^m\mid\sum_{k=0}^{\infty}g^*
 w_k(z)^* w_k(z)g <\infty\}. \label{d37}
\end{align}
Assume now that for some $z \in \BC_M$ there is a value $\wt \vp(z)$
such that 
\begin{equation} \label{d38}
\sum_{k=0}^{\infty}\begin{bmatrix} \wt  \vp(z)^* & I_{m_2}\end{bmatrix}w_k(z)^* w_k(z)
\left[
\begin{array}{c}
 \wt \vp(z) \\ I_{m_2}
\end{array}
\right] < \infty, \quad \wt \vp(z)\not=\vp(z).
\end{equation}
Then \eqref{d1!} and \eqref{d38} contradict the inequality \eqref{d37}.
Hence, there is no such $\wt \vp$ and the Weyl function is unique.
\epr

\paragraph{Weyl function: the inverse problem.} From Theorem \ref{TmD5} we know that the Weyl function of  a discrete skew-selfadjoint Dirac system \eqref{d1} with a pseudo-exponential potential generated  by a strongly admissible quadruple is a strictly proper rational matrix function. The next theorem shows  that the converse is also true.

\begin{thm} \label{TmD6}
Let $ \va $ be a strictly proper rational $m_1 \times m_2$ matrix function. Then $\va$ is the Weyl function of  a discrete skew-selfadjoint Dirac system \eqref{d1} with a pseudo-exponential  {potential} generated  by a strongly admissible quadruple $\si_0$. 
{A corresponding $ \si_0= \{\a, S_0,  \vt_1,   \vt_2\}$ can be obtained explicitly by using the  following procedure.
First, we construct a minimal realization  of $ \va:$
\begin{equation}\label{Dd1} 
\va(z)=-i \th_1^*(z I_n
+\g)^{-1}
\th_2.
\end{equation}
Next, we choose $X$ to be the unique positive definite solution of the Riccati equation
$\g X-X\g^*-iX\th_1\th_1^*X +i \th_2\th_2^*=0$. Finally, we put
\begin{equation} \label{Dd2}
S_0=I_n, \quad \vt_{1}=X^{ {1}/{2}} \theta_{1}, \quad \vt_{2}=X^{-{1}/{2}}  \theta_{2}, 
\quad  \a=X^{- {1}/{2}}\g  X^{ {1}/{2}}+i \vt_2\vt_2^*. 
\end{equation} }
\end{thm}
{ \bpr  Given a  strictly proper rational  matrix function $\va$, Corollary \ref{CyA5} shows that the procedure to recover $ \si_0= \{\a, S_0,  \vt_1,   \vt_2\}$ is well-defined
 and $\si_0$ is  strongly admissible. Moreover, according to Corollary \ref{CyA5}, $\vp$ admits representation
 \eqref{d24} where the quadruple $\{\a, S_0,  \vt_1,   \vt_2\}$ is given by
\eqref{Dd2}. Then we know from Theorem \ref{TmD5} that $\va$ is the Weyl function of system \eqref{d1}  with the potential $\{C_k\}$
generated by $\{\a, S_0,  \vt_1,   \vt_2\}$.}
\epr

\medskip

We conclude this section with some auxiliary results {on the $m\ts m$ matrices $\sH_k^+$ and $\sH_k^-$:
\begin{align}
 \sH_k^+ &:=2W_{\si_{k}}(i){P_1} W_{\si_{k}}(-i)^*, \quad {P_1}=(I_m+j)/2; \label{defHp1}\\
\sH_k^-&:=2W_{\si_{k}}(-i){P_2} W_{\si_{k}}(i)^*, \quad {P_2}=(I_m-j)/2, \label{defHm1} 
\end{align}
which will be essential in the next section. Here $W_{\si_{k}}$ is the transfer function of  the form \eqref{d7}.}
\begin{lem}\label{lemHpm1} Let $\si_0=\{ \a, S_0,  \vt_1, \vt_2\}$ be a strongly  admissible quadruple,   {let $\si_k$ be defined by \eqref{defsik}},
assume that $i\not \in \s(\a)$,  and let $\{C_k\}_{k\in \BN_0}$ 
be  the  pseudo-exponential potential generated by  $\si_0$. Then
\begin{align} 
\left(I_m+C_k \right )\sH_k^{-} &=\sH_{k+1}^{-} \left(I_m+C_k \right)=0, \label{m2a1}\\[.2cm]
 \left(I_m-C_k \right)\sH_k^{+} &=\sH_{k+1}^{+}\left(I_m-C_k \right)=0. \label{m2b1}
\end{align}
\end{lem}
\bpr  Since $\si_0=\{ \a, S_0,  \vt_1, \vt_2\}$ is admissible, we know (see Lemma \ref{lem11}) that $-i\not \in \s(\a)$. By assumption  $i\not \in \s(\a)$. This allows us to apply  \eqref{d19}, first  with $\lambda=i$ and next with $\lambda=-1$. It follows that both $W_{\si_{k}}(i)$ and $W_{\si_{k}}(-i)$ are well-defined and invertible. Moreover, we have
\begin{equation}
\label{invpm}
W_{\si_{k}}(i)^{-1}= W_{\si_{k}}(-i)^* \ands  W_{\si_{k}}(-i)^{-1}= W_{\si_{k}}(i)^*.
\end{equation}
Using these identities the formulas for $\sH_k^+$ and $\sH_k^-$ can be rewritten as:
\begin{equation}
\label{eqHpm1}
\sH_k^+=2W_{\si_{k}}(i)P_1W_{\si_{k}}(i)^{-1}, \quad  \sH_k^-=2W_{\si_{k}}(-i)P_2W_{\si_{k}}(-i)^{-1}.
\end{equation}

Next, applying \eqref{d9}, first with $\lambda=i$ and next with $\lambda =-i$,   we see that
\begin{align} 
&(I_m-C_k)W_{\si_{k}}(i)=2W_{\si_{k+1}}(i){P_2}, \label{d9p}\\[.2cm]
&(I_m+C_k)W_{\si_{k}}(-i)=2W_{\si_{k+1}}(-i){P_1}.  \label{d9m}
\end{align}
It follows that
\begin{align} \nn
(I_m-C_k)\sH_k^+&=2(I_m-C_k)W_{\si_{k}}(i){P_1} W_{\si_{k}}(-i)^*\\[.2cm]
\label{v1}
&= 4W_{\si_{k+1}}(i){P_2}{P_1} W_{\si_{k}}(-i)^*=0, \\[.2cm] \nn
(I_m+C_k)\sH_k^-&=2(I_m+C_k)W_{\si_{k}}(-i){P_2} W_{\si_{k}}(i)^*\\[.2cm]
&= 4W_{\si_{k+1}}(-i){P_1}{P_2} W_{\si_{k}}(i)^*=0. \label{v2}
\end{align} 
From \eqref{eqHpm1}, we obtain 
\begin{align}
&\sH_k^+W_{\si_{ k}}(i)=2W_{\si_{k}}(i){P_1} W_{\si_{k}}(i)^{-1}W_{\si_{ k}}(i)=2W_{\si_{k}}(i){P_1},\label{m42p}\\[.2cm]
&\sH_k^-W_{\si_{ k}}(-i)=2W_{\si_{k}}(-i){P_2} W_{\si_{k}}(-i)^{-1}W_{\si_{ k}}(-i)=2W_{\si_{k}}(-i){P_2}.\label{m42m}
\end{align} 
In particular, we have 
\begin{align} \label{v3}
\sH_k^+W_{\si_{k}}(i){P_2}=0 \ands \sH_k^-W_{\si_{k}}(-i){P_1}=0 \quad (k \geq 0).
\end{align} 
Using \eqref{d9p} and \eqref{d9m} again {and taking into account \eqref{v3}},  we see that
\begin{align} \nn
\sH_{k+1}^+({I_m}-C_k)W_{\si_{ k}}(i)&=\sH_{k+1}^+({I_m}-C_k)W_{\si_{ k}}(i)({P_1}+{P_2})\\ \nn
&=2\sH_{k+1}^+W_{\si_{ k+1}}(i){P_2}({P_1}+{P_2})\\ 
&=2\sH_{k+1}^+W_{\si_{ k+1}}(i){P_2}={0},  \label{v4}
\end{align}
and 
\begin{align} \nn
\sH_{k+1}^-({I_m}+C_k)W_{\si_{ k}}(-i)&=\sH_{k+1}^-({I_m}+C_k)W_{\si_{- k}}(-i)({P_1}+{P_2})\\ \nn
&=2\sH_{k+1}^-W_{\si_{ k+1}}(-i){P_1}({P_2}+{P_1})\\
&=2\sH_{k+1}^-W_{\si_{ k+1}}(-i){P_1}={0}.  \label{v5}
\end{align}
{Formulas  \eqref{v1}, \eqref{v2}, \eqref{v4} and \eqref{v5} yield} \eqref{m2a1} and \eqref{m2b1}.\epr

\medskip
For later purposes (see the next section) we mention that under the assumptions of Lemma \ref{lemHpm1} the matrices $\sH_k^+$ and $\sH_k^-$ can be rewritten as
\begin{align} 
& \sH_k^+ = I_m+W_{\si_{ k}}(i)jW_{\si_{ k}}(-i)^*,\label{defHp2} \\[.1cm]
&\sH_k^- =I_m-W_{\si_{ k}}(-i)jW_{\si_{ k}}(i)^*. \label{defHm2}
\end{align}
Indeed, these identities follow from \eqref{defHp1} and \eqref{defHm1} by using the two indenties in \eqref{invpm}. Together \eqref{defHp2} and \eqref{defHm2} imply that
\[
\sH_k^+ + (\sH_k^-)^*=2I_m.
\]

Finally, using \eqref{defHp2} we see that
\begin{align}
\sH_k^+ &= I_m+\left(I_m+i\la_k^*S_k^{-1}(iI_n-\a)^{-1}\la_k\right)j \nn \\
& \hspace{2cm}\ts \left(I_m+i\la_k^*S_k^{-1}(-i I_n-\a)^{-1}\la_k\right)^* \nn \\
&= I_m+\left(j+i\la_k^*S_k^{-1}(iI_n-\a)^{-1}\la_k j\right) \nn \\
& \hspace{2cm}\ts \left(I_m -i\la_k^*(i I_n-\a^*)^{-1}S_k^{-1}\la_k\right ) \nn \\
&= I_m+j+i\la_k^*S_k^{-1}(iI_n-\a)^{-1}\la_kj\nn\\
& \hspace{2cm}- ij \la_k^*(i I_n-\a^*)^{-1}S_k^{-1}\la_k\nn\\
& \hspace{2cm}+\la_k^*S_k^{-1}(iI_n-\a)^{-1}\la_kj\la_k^*(i I_n-\a^*)^{-1}S_k^{-1}\la_k.\label{Hkpl2}
\end{align}
Analogously, using \eqref{defHm2}, we have 
\begin{align}
{\sH_k^-} &= I_m-\left(I_m+i\la_k^*S_k^{-1}(-iI_n-\a)^{-1}\la_k\right)j \nn \\
& \hspace{2cm}\ts \left(I_m+i\la_k^*S_k^{-1}(i I_n-\a)^{-1}\la_k\right)^* \nn \\ 
&= I_m-\left(j+i\la_k^*S_k^{-1}(-iI_n-\a)^{-1}\la_k j\right) \nn \\
& \hspace{2cm}\ts \left(I_m -i\la_k^*(-i I_n-\a^*)^{-1}S_k^{-1}\la_k\right ) \nn \\
&= I_m- j+i\la_k^*S_k^{-1}(iI_n+\a)^{-1}\la_k j\nn \\
& \hspace{2cm}- ij \la_k^*(i I_n+\a^*)^{-1}S_k^{-1}\la_k\nn\\
& \hspace{2cm}-\la_k^*S_k^{-1}(iI_n+\a)^{-1}\la_kj\la_k^*(i I_n+\a^*)^{-1}S_k^{-1}\la_k.\label{Hkml2}
\end{align}

\section{Generalized discrete Heisenberg magnet \\ model {and its auxiliary linear systems}} \label{NL}
\setcounter{equation}{0}

{This section is devoted to the generalized discrete Heisenberg magnet model.} {We apply the quadruples considered in the previous section in order
to construct explicit solutions of this model.  Given potential $\{C_k(t)\}$ from such solution and using the results from the previous section, we  also express explicitly evolution of the Weyl function
of the auxiliary linear system, which coincides with system \eqref{d1}.}
\subsection{Generalized discrete Heisenberg magnet  model}
For the case that $m_1=m_2=1$,  system \eqref{d1} is an auxiliary linear system for the integrable isotropic Heisenberg magnet model \cite{Skl} (see also \cite[Part II, Section 1.2]{FT}).  
More precisely, it was shown in \cite{FT, Skl} that   the isotropic Heisenberg magnet model is equivalent to the compatibility condition
\begin{align}\label{m4}&
\frac{d}{dt}G_k=F_{k+1}G_k-G_kF_k
\end{align}
of the auxiliary systems
\begin{align}\label{m5}&
y_{k+1}=G_k y_k, \quad G_k(t,z):=I_m+\frac{i}{z}C_k(t);
\\ \label{m6}& \frac{d}{dt} y_k=F_k y_k, \quad F_k(t,z):=-\frac{H_k^{+}(t)}{z+i} -\frac{H_k^{-}(t)}{z-i},
\end{align}
where $m=2$,
\begin{align}
& H_k^{\pm}(t)=f_k(t)(I_2 \pm C_k(t))(I_2\pm C_{k-1}(t)), \label{mg2}
\end{align}
$f_k$ is a scalar function depending on $C_{k-1}$ and $C_k$, $k>0$, $C_k$ has the form
\begin{align}\label{vm1}&
C_k=\begin{bmatrix}c_k^3 &
c_k^1-ic_k^2
\\ c_k^1+i c_k^2 & -c_k^3
\end{bmatrix}, \quad (c_k^1)^2+(c_k^2)^2+(c_k^2)^2=1,
\end{align} 
and the numbers  $c_k^1$, $c_k^2$ and  $c_k^3$ are real valued. It is easy to see that representation \eqref{vm1} is equivalent to the conditions
$C_k=C_k^*$, $\det C_k=1$ and Tr~$C_k=0$, where Tr means trace. These conditions are in turn equivalent to the conditions $C_k=C_k^*$
and the eigenvalues of $C_k$ equal $1$ and $-1$. Hence, \eqref{vm1} is equivalent to the second part of \eqref{d1}, where $m_1=m_2=1$ (see \cite{KaS}).
Moreover, the second part of \eqref{d1}  yields relations
\begin{align}\label{vm2}&
I_2+C_k=2U_kP_1U_k^*, \quad I_2-C_k=2U_kP_2U_k^*, 
\end{align} 
where $P_1=(I_2+j)/2$,  $P_2=(I_2-j)/2$ and $P_1P_2=P_2P_1=0$.  Therefore, \eqref{mg2} and \eqref{vm2} imply that $H_k^{\pm}$ has rank 1, that
\begin{align}\label{vm3}
 &
(I_2 \mp C_k(t))H_k^{\pm}(t)=H_k^{\pm}(t)(I_2 \mp C_{k-1}(t))=0,
\end{align} 
and that \eqref{vm3} defines matrices $H_k^{+}$ and $H_k^{-}$ up to scalar factors.

Generalized Heisenberg magnet models are also actively studied
(see, e.g., \cite{DimM, Gol, KuR, SalH, Skr, Suth}), including the case of  the continuous generalized Heisenberg magnet model, where the $m\times m$ matrix function $C$ (an analog of $C_k$) satisfies the condition $C=C^*=C^{-1}$ (see \cite[Remark 5]{DimM} and references therein). Using an analog of condition \eqref{vm3}
instead of \eqref{mg2}, we consider the following  generalized discrete Heisenberg magnet model $($DGHM$)$:
\begin{align} 
&i\frac{d}{dt}C_k(t) =\left(H_{k+1}^{-}(t)-H_{k+1}^{+}(t)\right)C_k(t)\nn \\
 &\hspace{5cm}-C_k(t)\left(H_k^{-}(t)-H_k^{+}(t)\right),\label{m1}
\end{align}
where 
\begin{align}
&\left(I_m+C_k(t)\right )H_k ^{-}(t)=H_{k+1}^{-}(t)\left(I_m+C_k(t)\right)=0, \label{m2a}\\
 &\left(I_m-C_k(t)\right)H_k ^{+}(t)=H_{k+1}^{+}(t)\left(I_m-C_k(t)\right)=0,\label{m2b}
\end{align}
and
\begin{align} 
 &  C_k(t)=U_k(t)jU_k(t)^*, \quad U_k(t) U_k(t)^*=I_m,\quad k\in \BN_0.\label{mg3a}
\end{align}
The next proposition shows that \eqref{m1} coincides with \eqref{m4}.
\begin{Pn}\label{mPnComp}  Under conditions \eqref{m2a}--\eqref{mg3a},
the discrete generalized  \\  
Heisenberg magnet model  \eqref{m1}
is equivalent to the compatibility condition \eqref{m4}
of the auxiliary systems \eqref{m5} and \eqref{m6}.
\end{Pn}
\bpr  In view of  the definitions of $G_k$ and $F_k$, we rewrite \eqref{m4} in the form
\begin{align}\label{m7}&
\frac{i}{z}\frac{d}{dt}C_k=\left(I_m+\frac{i}{z}C_k\right)\left(\frac{H_k^{+}}{z+i} +\frac{H_k^{-}}{z-i}\right)-
\left(\frac{H_{k+1}^{+}}{z+i} +\frac{H_{k+1}^{-}}{z-i}\right)\left(I_m+\frac{i}{z}C_k\right).
\end{align}
Using the equality $\frac{\pm i}{z(z\pm i)}=\frac{1}{z}-\frac{1}{z\pm i}$, after simple transformations we see that \eqref{m7} is equivalent to the relation
\begin{align}\nn
i\frac{d}{dt}C_k=& (H_{k+1}^{-}-H_{k+1}^{+})C_k-C_k(H_k^{-}-H_k^{+})
\\  \nn \,\, &+\frac{z}{z+i}\big((I_m-C_k)H_k^{+}-H_{k+1}^{+}(I_m-C_k)\big)
\\ \label{m8} \,\, &
+\frac{z}{z-i}\big((I_m+C_k)H_k^{-}-H_{k+1}^{-}(I_m+C_k)\big).
\end{align}
Taking into account the equalities \eqref{m2a} and \eqref{mg3a} ,
we derive that  the terms in the second and third lines of \eqref{m8} both turn to zero,
and therefore \eqref{m8} is equivalent to \eqref{m1}.
\epr

In order to construct explicit solutions of DGHM, we introduce an additional parameter $t$
into the quadruples generating potentials $\{C_k(t)\}$. Namely, we consider  quadruples $\si_t=\si_{t,0}= \{\a, \, S_0(t), \, \vt_1(t), \vt_2(t)\}$.
Here $\a$ is constant and the dependence of the quadruples on $t$ is described by the equations
\begin{align}\label{m9!}
\left(\frac{d}{dt}\Lam_0\right)(t)=&-2\big((\a-i I_n)^{-1}\Lam_0(t)P_1 + (\a+i I_n)^{-1}\Lam_0(t)P_2\big), 
\\ 
\label{m9'} P_{1}:=&(I_m + j)/2, \quad P_{2}:=(I_m - j)/2;
\\
\nn 
\left(\frac{d}{dt}S_0\right)(t)=&
- \Big( (\a -i I_n)^{-1} S_0(t)+(\a
+i I_n)^{-1} S_0(t)
\\  \nonumber &
\quad + S_0(t)(\a^* +i I_n)^{-1}
 + S_0(t)(\a^*
-i I_n)^{-1}
\\
\nonumber & \quad 
+ 2(\a^2+I_n)^{-1} 
\big( \a \Lam_0(t)j
\Lam_0(t)^*+\Lam_0(t)j \Lam_0(t)^* \a^* \big)
\\  & \qquad \qquad  \qquad  \qquad  \qquad  \qquad  \qquad  \,\,\,
\times ( (\a^*)^2+I_n)^{-1} \Big). \label{m10}
\end{align}
Formula \eqref{m9!} yields the equalities
\begin{align}  
&\hspace{2.2cm}\la_0(t)=\begin{bmatrix} \vt_1(t) &\vt_2(t) \end{bmatrix}, \label{m9}\\
&\vt_1(t)=e^{-2t(\a-{i} I_n)^{-1}}\vt_1(0),      \quad  \vt_2(t) =e^{-2t(\a+{i} I_n)^{-1}}\vt_2(0).\label{m15}
\end{align}
We assume for simplicity that $\si_{0,0}$ is strongly admissible and that $i \not\in \s(\a)$.
Then, the quadruples $\si_{t,0}$ are well-defined and  uniquely determined by the initial quadruple  $\si_{0,0}$.
Moreover, all $\si_{t,0}$ are strongly admissible (see Subsection \ref{subsec4.2}).
It was shown in Section \ref{ddies} that   the strongly admissible quadruples  $\si_{t,0}$
determine  strongly admissible quadruples  $\si_{t,k}$
using relations \eqref{defLak} and \eqref{defSk},  
and that  in this way $\si_{t,0}$ generate potentials $\{C_k(t)\}$.
Since $\si_{t,0}$ are determined by $\si_{0,0}$ we also say that $\si_{t,k}$ are determined by $\si_{0,0}$ and that potentials $\{C_k(t)\}$ are generated by  $\si_{0,0}$.
The following theorem is our main result in this section.
\begin{Tm} \label{thm2var}  Let   $\si_{0,0}=\{ \a, S_0,  \vt_1, \vt_2\}$   be a strongly  admissible quadruple, and assume that $  i\not \in \s(\a)$.  Then for all $t\in \BR$   and $k\in \BN_0$   the quadruples  $\si_{t, k}$ determined by $\si_{0,0}$ are strongly admissible.  Furthermore, 
the matrix functions
\begin{align}\label{defHplus}
H_k^+(t):=&\sH_k^+(t)= W_{\si_{t,k}}(i)(I_m+j) W_{\si_{t,k}}(-i)^*,
\\  \label{defHmin}
H_k^-(t):=&\sH_k^-(t)=W_{\si_{t,k}}(-i)(I_m-j) W_{\si_{t,k}}(i)^*
\end{align}
are well-defined and satisfy DGHM \eqref{m1} and conditions \eqref{m2a} and \eqref{m2b}.
\end{Tm}
Recall that matrix functions $\sH_k^{\pm}$ were introduced and studied in the last part of Section \ref{ddies}.
The proof of Theorem \ref{thm2var} and some related results on discrete Dirac systems depending on an additional continuous time parameter
are given in the next subsection.

\medskip\noindent
Finally, we note that in view of \eqref{m15} we have the following corollary of Theorems \ref{TmD5} and \ref{thm2var}.
\begin{cor} \label{CyDGHM} Under the conditions of Theorem
\ref{thm2var}, the evolution of the Weyl function $\varphi(t,z)$ of
the system \eqref{m5} is given by the
formula
\begin{align} \nn &
\varphi (t, z)=-{i} \vt_1^* {e}^{-2t(\a^* +{i} I_n)^{-1}}S_0(t)^{-1}(z
I_n + \b(t))^{-1} e^{-2t(\a +{i} I_n)^{-1}} \vt_2,
\\ \nn &
\b(t): =\a -{i} {e}^{-2t(\a +{i} I_n)^{-1}} \vt_2 \vt_2^* e^{-2t(\a^*
-{i} I_n)^{-1}} S_0(t)^{-1}, \,\, \vt_k:=\vt_k(0) \,\, (k=1,2).
\end{align}
\end{cor}

\subsection{Discrete {Dirac} systems depending on {an additional continuous time} parameter}\label{subsec4.2}
Let   $\si=\{ \a, S_0,  \vt_1, \vt_2\}$   be a strongly admissible quadruple, and assume that $ i\not \in \s(\a)$.   We shall associate  with $\si$  a family of  strongly admissible quadruples depending on two parameters $t$ and $k$, where $t\in \BR_+$ and $k\in \BN_0$. This will be done in two steps. 

\medskip
\noindent\textsc{Step 1: construction of $\si_t$.} In this first step we just  require $\si$ to be admissible, not necessarily strongly admissible.  From item (i) in Lemma \ref{lem11} we know that  $- i\not \in \s(\a)$. Together with the assumption that  $ i\not \in \s(\a)$ this implies that both  $\a-{i} I_n$ and $\a+{i} I_n$ are invertible.  {This allows us to define $\Lam(t)=\Lam_0(t)$ and $S(t)=S_0(t)$
using \eqref{m15} (and \eqref{m9}) and \eqref{m10}, respectively.}

\begin{lem} \label{lemdiscr1a} Let   $\si=\{ \a, S_0,  \vt_1, \vt_2\}$   be an admissible quadruple, and assume that $  i\not \in \s(\a)$.  Let $\vt_j(t) $, $j=1,2$, and $S(t)$ be defined by \eqref{m15} and \eqref{m10}, respectively. Then on some interval  interval $-\ve_1<t<\ve_2$ the quadruple  
\begin{equation}\label{Sigma[t]} 
\si_{t}=\{ \a, S(t),  \vt_1(t), \vt_2(t)\}
\end{equation}
is   admissible too. Moreover, if $\si$ is strongly admissible,  then the same holds true for the quadruple $\si_{t}$ {for each $t\in \BR$}.
\end{lem}
\bpr Since $S_0$ is positive definite, the same holds true for $S(t)$  provided $t$ is sufficiently small. Put
\begin{equation}\label{m11} 
\Up(t):=\a S(t)-S(t)\a^*-i\Lam(t)\Lam(t)^*.
\end{equation}
 Here $\Lam(t)$ is defined by \eqref{m9}. Since $\{ \a, S_0,  \vt_1, \vt_2\}$ is an admissible  quadruple, it follows that $\Up(0)=0$.  After some easy transformations, using \eqref{m9} and \eqref{m10}  and differentiating the right-hand side of \eqref{m11}, we see that $ \Up $ satisfies  the linear differential equation
\begin{align} 
\frac{d}{d t}\Up(t)&=-(\a+i I_n)^{-1}\big(\a \Up(t)+\Up(t)\a^*\big)(\a^*-i I_n)^{-1}+ \nn\\  
&\hspace{2cm}-(\a-i I_n)^{-1}\big(\a \Up(t)+\Up(t)\a^*\big)(\a^*+i  I_n)^{-1}.\label{m12}
\end{align}
Since the initial value $\Up(0)$ is zero, it follows  that $\Up(t)\equiv 0$, that is, 
\begin{equation}\label{m13} 
\a S(t)-S(t)\a^*=i\Lam(t)\Lam(t)^*.
\end{equation}
Using  \eqref{m15} and \eqref{m9}, we conclude that $\{ \a, S(t),  \vt_1(t), \vt_2(t)\}$ is an admissible quadruple. 

Next, assume $\si$   is strongly  admissible.  {The identity \eqref{m13} was already proved above.}  
By Definition \ref{DnPE},  the fact that $\si$ is strongly admissible means that the pair $\a$ and $\vt_1$ is controllable.  But then the first equality in \eqref{m15} tells us that
\begin{align}\label{m16}&
\spa \bigcup_{\nu=0}^{n-1}\im \a^\nu \vt_1(t)=e^{-2t(\a-i I_n)^{-1}}\Big(\spa \bigcup_{\nu=0}^{n-1}\im\big(\a^\nu \vt_1(0) \Big)=\BC^n,
\end{align}
that is,  the pair $\a$ and $\vt_1(t)$ is also controllable. {Furthermore, since $\si$ is strongly admissible item (ii)
of Lemma \ref{lem11} shows that $\s(\a)\subset \BC_+$. From \eqref{m13} and $\s(\a)\subset \BC_+$ 
follows representation \eqref{altS0} of $S_0(t)=S(t)$. Since the pair $\a$ and $\vt_1(t)$ is controllable
representation \eqref{altS0} implies that $S(t)>0$.}
Hence, again using  Definition~\ref{DnPE}, the quadruple $\si_{t}$ is strongly admissible. 
\epr

\medskip
\noindent\textsc{Step 2: construction of $\si_{t,k}$.} 
Let  $\si=\{ \a, S,  \vt_1, \vt_2\}$   be  a strongly admissible quadruple, and assume that $  i\not \in \s(\a)$.  Let $\si_t$ be the quadrupple defined by \eqref{Sigma[t]}. {Hence}, according to the previous lemma,  the quadruple  $\si_{t}$ is also  strongly admissible for each {$t\in \BR$.}   But then we can apply the results of Section \ref{ddies} to associate with each $\si_t$ a family of 
quadruples  $\si_{t, k}$, $k\in \BN_0$. 

Since $\si_{t}$ is strongly admissible, we know from Lemma \ref{lemadmdiscrk} that $\a$ is invertible. In fact, by   item (ii) in Lemma~\ref{lem11},  all eigenvalues of $\a$ are in $\BC_+$. Following the constructions given in Section~\ref{ddies}, put
\begin{align}
&\vt_{1}(t, k)=(I_n+i\a^{-1})^k\vt_1(t), \quad  \vt_{2, k}(t)=(I_n-i\a^{-1})^k\vt_{2}(t), \label{m15a}\\[.2cm]
&\hspace{2cm}\la(t, k)= \begin{bmatrix}\vt_{1}(t, k) &\vt_{2}(t, k) \end{bmatrix}.\label{m9a}
\end{align}
Furthermore, define $S(t,k)$ by 
\begin{align*}
&S(t, k+1)=S(t,k) +\a^{-1}S(t,k) \a^{-*}+ \\
&\hspace{3cm}+\a^{-1}S(t,k) \a^{-*} \la(t, k) j \la(t, k)^*,  \quad k=0,1,2, \ldots,\\
&S(t, 0)=S(t).
\end{align*}
Using these matrices we set 
\begin{equation}\label{defSigtk}
\si_{t, k}=\{ \a, S(t,k),  (I_n+i\a^{-1})^k\vt_1(t), (I_n-i\a^{-1})^k\vt_2(t)\}.
\end{equation}
By Lemma \ref{lemadmdiscrk}, the quadruple $\si_{t, k}$ is also strongly admissible.

\medskip
Let  $\si=\{ \a, S,  \vt_1, \vt_2\}$   be  a strongly admissible quadruple, and  let $\si_t$ be the strongly admissible quadruple defined by \eqref{Sigma[t]}.  Following \eqref{d5} put 
\begin{align*}
C_k(t)& =j+\la(t,k)^*S(t,k)^{-1}\la(t,k)+\\
&\hspace{2cm} - \la(t,k+1)^*S(t,k+1)^{-1}\la(t,k+1), \quad k\in \BN_0.
\end{align*}
Then, by Definition \ref{DnPE},  the sequence $\{C_k(t)\}_{k\in \BN_0}$  is the pseudo-expontential potential generated by $\si_t$. The corresponding  SkDDS is given by 
\[
y_{k+1}(t, z)=\Big(I_m+iz^{-1}C_k(t)\Big)y_{k}(t, z) \quad (k\in \BN_0).
\]
Furthermore, we know (see Propositon \ref{PnDa1}) that $ C_k(t)$ admits a factorization of the form:
\begin{equation}\label{m3}
 C_k(t)=U_k(t)jU_k(t)^*, \quad U_k(t) U_k(t)^*=I_m, \quad k\in \BN_0.
\end{equation}
{We also  note} that \eqref{m2a} and \eqref{m2b} are immediate consequences  of  the identities  \eqref{m2a1} and \eqref{m2b1}. Indeed, since $\si_{t}=\{ \a, S (t),  \vt_1(t), \vt_2(t)\}$  is   strongly admissible, we just apply Lemma \ref{lemHpm1} with $\si_{t}=\si_{t,0}$ in place of $\si_0$. Thus it remains to prove the identity \eqref{m1} {in order to prove Theorem \ref{thm2var}}.
{First, we formulate a proposition, which will be proved later.}
\begin{prop}\label{propComp2}  Let   $\si=\{ \a, S_0,  \vt_1, \vt_2\}$   be a strongly  admissible quadruple, assume that $  i\not \in \s(\a)$, and let  $\si_{t, k}$ be the strongly admissible quadruple defined by \eqref{defSigtk}, where $k\in \BN_0$ and {$t\in \BR$.} Then the function  $Y_k(t,z)$   defined by
\begin{equation}   \label{m44} 
Y_k(t,z)=W_{\si_{t, k}}(-z)\left(I_m+\frac{{i}}{z}j\right)^k\exp\left\{-2t\left(\frac{P_1}{z+{i}}+\frac{P_2}{z-{i}}\right)\right\}, 
\end{equation}
where $z\not=\pm {i}$ and $z\not\in\s(-\a)$, is   a solution of   equation \eqref{m5} as well as of  equation \eqref{m6}.
\end{prop}\medskip
\noindent\textbf{Proof of Theorem \ref{thm2var}.} As soon as Proposition \ref{propComp2}  is proved we can prove Theorem \ref{thm2var}  by using  Proposition \ref{mPnComp}. Indeed, assume that the function $Y$ defined by \eqref{m44} satisfies \eqref{m5} and \eqref{m6}. Then
\begin{equation*}
\frac{d}{dt} Y_{k+1}(t,z)  =F_{k+1}(t,z)Y_{k+1}(t,z)=F_{k+1}(t,z)G_k(t,z)Y_{k}(t,z),
\end{equation*}
and 
\begin{align*}
\frac{d}{dt} Y_{k+1}(t,z)& =\frac{d}{dt}\Big(G_k(t,z)Y_{k}(t,z)\Big)\\[.2cm]
&=\left(\frac{d}{dt}G_k(t,z)\right)Y_{k}(t,z)+G_k(t,z)\frac{d}{dt} Y_{k}(t,z)\\[.2cm]
&=\left(\left(\frac{d}{dt}G_k(t,z)\right)+G_k(t,z)F_k(t,z)\right)Y_{k}(t,z).
\end{align*}
We conclude that
\[
 \left(\frac{d}{dt}G_k(t,z)\right)Y_{k}(t,z)=\Big(F_{k+1}(t,z)G_k(t,z)-G_k(t,z)F_k(t,z)\Big)Y_{k}(t,z).
\]
Finally, note that $z\not=\pm {i}$ and $z\not\in\s(-\a)$ imply that   $Y_{k}(t,z)$ is non-singular.  But then the above identity yields  the compatibility equation \eqref{m4}, and we can apply Proposition  \ref{mPnComp} to get the identity \eqref{m1}.
\epr

\medskip
In order to prove Proposition \ref{propComp2} we need to  compute the derivatives
\begin{equation}
\label{dervas}
\frac{d}{d t} \la(t,k), \quad \frac{d}{d t}S(t,k), \quad \frac{d}{d t}  \la(t,k)^* S(t,k)^{-1}, \quad   \frac{d}{d t}W_{\si_{t, k}}.
\end{equation}
This will be done in a couple of steps. 

\medskip
\noindent\textsc{The first derivative in \eqref{dervas}.} Recall that $\la(t,k)$ is given by \eqref{m15a} and \eqref{m9a}.  It follows that
\[
\frac{d}{d t} \la(t,k)=
\begin{bmatrix} (I_n+i\a^{-1})^k\big(\frac{d}{d t}\vt_1(t)\big)&(I_n-i\a^{-1})^k\big(\frac{d}{d t}\vt_2(t)\big)\end{bmatrix}.
\]
But then, using \eqref{m15} and \eqref{m9}, we see that 
\begin{align}
&\frac{d}{d t} \la(t,k)=\nonumber\\[.1cm]
&\hspace{.7cm}=-2\begin{bmatrix} (\a-i I_n)^{-1} (I_n+i\a^{-1})^k\vt_1(t)& (\a+i I_n)^{-1}(I_n-i\a^{-1})^k\vt_2(t)\end{bmatrix}\nonumber\\[.1cm]
&\hspace{.7cm}= -2\begin{bmatrix}(\a-i I_n)^{-1} \vt_{1}(t,k)& (\a+i I_n)^{-1}  \vt_{2} (t,k)\end{bmatrix}\nonumber\\[.1cm]
&\hspace{.7cm}=-2 \Big( (\a -i I_n)^{-1} \la(t,k) {P_1}  + (\a
+i I_n)^{-1} \la(t,k)  {P_2} \Big).\label{m18}
\end{align}
Here ${P_1}$ and ${P_2}$ are the projctions defined in \eqref{defHplus} and  \eqref{defHmin}, respectively.

\medskip
\noindent\textsc{The second derivative in \eqref{dervas}.}  In order to compute the second derivative in  \eqref{dervas} we present an alternative way to obtain the quadruple $\si_{t,k}$ given by \eqref{defSigtk}. As before, let    $\si=\{ \a, S_0,  \vt_1, \vt_2\}$   be a strongly admissible quadruple, and assume that  $i\not \in \s(\a)$.  We know that $\s(\a)\subset \BC_+ $. In a particular $\a$ is non-singular.  As in   Section \ref{ddies},  let
\[
\si_k= \{\a, S_k, (I_n+i \a^{-1})^{k}\vt_1,  (I_n-i \a^{-1})^{k}\vt_2\}.
\]
From Lemma \ref{lemadmdiscrk}  we know that $\si_k$ is again a strongly admissible quadruple. The assumption $i\not \in \s(\a)$ allows us to apply Lemma \ref{lemdiscr1a} to  $\si_k$ in place of $\si$.   Put 
\begin{align}
\tilde{\vt}_1(k, t)&= e^{-2t(\a-{i} I_n)^{-1}} (I_n+i \a^{-1})^{k}\vt_1, \label{vt1kt}\\[.2cm]
\tilde{\vt}_2(k, t)&=e^{-2t(\a+{i} I_n)^{-1}} (I_n-i \a^{-1})^{k}\vt_1, \label{vt2kt}\\[.2cm]
\tilde{\Lam}(k, t)&= \begin{bmatrix} \tilde{\vt}_1(k, t)&  \tilde{\vt}_2(k, t)
\end{bmatrix}, \label{tlakt}
\end{align}
and let $\tilde{S}(k, t)$ be given by 
\begin{align}
& \frac{d}{dt}\tilde{S}(k, t) =\nn\\ 
&=  - \Big[ (\a -{i} I_n)^{-1}\tilde{S}(k, t)+(\a+{i} I_n)^{-1} \tilde{S}(k, t)+ 
\nn\\ 
& \quad \hspace{2.5cm} + \tilde{S}(k, t)(\a^* +{i} I_n)^{-1}+\tilde{S}(k, t)(\a^* -{i} I_n)^{-1}+ \label{m19a}  \\
& \quad  + 2(\a^2+I_n)^{-1}\Big( \a \tilde{\Lam}(k,t)j
\tilde{\Lam}(k,t)^*+\tilde{\Lam}(k,t)j \tilde{\Lam}(k,t)^* \a^* \Big)\left( (\a^*)^2+I_n\right)^{-1} \Big],\nn \\
&\tilde{S}(k, 0)=S_k. \nn 
\end{align}
Note that \eqref{m19a} is the analogue of \eqref{m10} with $\si_k$ in place of $\si$. Thus, by Lemma \ref{lemdiscr1a}, the  quadruple
\[
\si_{k,t} = \{ \a, \tilde{S}(k, t),  \tilde{\vt}_1(k, t),  \tilde{\vt}_2(k, t)\}  
\]
is strongly admissible.

\begin{lem}\label{lemcompat} For each $k\in \BN_0$ and $t \in \BR$  we have  $\si_{t, k}= \si_{k, t} $. In particular,
\begin{align}  
& \frac{d}{dt}S(t,k)=\nn\\
&=  - \Big[ (\a -{i} I_n)^{-1} S(t,k)+(\a+{i} I_n)^{-1} S(t,k)+\nn \\ 
& \quad \hspace{2.5cm} +  S (t,k)(\a^* +{i} I_n)^{-1}+S(t,k)(\a^* -{i} I_n)^{-1}+  \label{m19} \\
& \quad  + 2(\a^2+I_n)^{-1}\Big( \a \la(t,k)j
\la(t,k)^*+\la(t,k)j \la(t,k)^* \a^* \Big)( (\a^*)^2+I_n)^{-1} \Big].\nn
\end{align}

\end{lem}
\bpr From \eqref{m15} and \eqref{m15a}  and the identities \eqref{vt1kt} and \eqref{vt2kt} we see that 
\begin{align*}
\vt_1(t,k)&=\tilde{\vt}_1(k, t) =e^{-2t(\a-{i} I_n)^{-1}}(I_n+i \a^{-1})^{k}\vt_1,\\
\vt_2(t,k) &= \tilde{\vt}_2(k, t)=e^{-2t(\a+{i} I_n)^{-1}} (I_n-i \a^{-1})^{k} \vt_2.
\end{align*}
It follows that $\Lam(t,k)=\tilde{\Lam}(k,t)$.  Recall that $\si_{t,k}$ and $\si_{k, t}$ are both strongly admissible. Then item (ii) in Lemma \ref{lem11} tells us that $S(t, k$ and $\tilde{S}(k, t)$ are uniquely determined by $\Lam(t,k)$ and $\tilde{\Lam}(k,t)$.  Thus $S(t, k$ and $\tilde{S}(k, t)$ coincide, and hence $\si_{t, k}= \si_{k, t} $. Thus $\si_{t, k}= \si_{k, t} $.  

Finally, using $\Lam(t,k)=\tilde{\Lam}(k,t)$ and $S(t, k)=\tilde{S}(k, t)$, we see that \eqref{m19} follows from \eqref{m19a}. \epr

 \medskip
\noindent\textsc{The third derivative in \eqref{dervas}.}  From \eqref{m18} and \eqref{m19} we shall derive the identity
\begin{align} 
&\frac{d}{dt} \left( \la(t,k)^* S(t,k)^{-1} \right)=H_k^+(t)  \la(t,k)^* S(t,k)^{-1}(\a - {i} I_n)^{-1} +\nn\\
&\hspace{5cm}+H_k^-(t)\la(t,k)^* S(t,k)^{-1}(\a +{i} I_n)^{-1}. \label{m28}
\end{align}
Here $H_k^+$ and $H_k^-$ are the functions given by \eqref{defHplus} and \eqref{defHmin}, respectively. To do this, note that  according to  (\ref{m18}) we have
\begin{align} \nonumber
\frac{d}{dt} \left( \la(t,k)^*  S(t,k)^{-1} \right) & =-  2{P_1}
\la(t,k)^*(\a^*+i I_n)^{-1} +\nn \\ 
&\qquad- 2{P_2} \la(t,k)^*(\a^*-{i} I_n)^{-1}+\nn\\
&\qquad  -\la(t,k)^*S(t,k)^{-1} \Big(\frac{ d }{dt}S(t,k) \Big) S(t,k)^{-1}.\label{m31}
\end{align}
Identity \eqref{d2'} implies that
\begin{align}  
&(\a^* \pm {i} I_n)^{-1}S(t,k)^{-1}=S(t,k)^{-1}(\a \pm {i}
I_n)^{-1}+\nn \\
&\hspace{1cm}+{i} (\a^* \pm {i} I_n)^{-1}S(t,k)^{-1} \la(t,k) \la(t,k)^* S(t,k)^{-1}(\a \pm {i} I_n)^{-1}.\label{m32}
\end{align}
By using (\ref{m19}), (\ref{m32}) and the  equalities
\begin{align} \nn &
2( \a^2 + I_n)^{-1}={i}\big((\a + {i} I_n)^{-1}- (\a -{i} I_n)^{-1}\big), 
\\ \label{m33}
&  2( \a^2 + I_n)^{-1} \a= (\a - {i} I_n)^{-1}+ (\a +{i} I_n)^{-1},
\end{align}
after some calculations, we see that  (\ref{m31}) can be rewritten  in the form:
\begin{align*} 
&\frac{d}{dt} \left( \la(t,k)^* S(t,k)^{-1} \right)= \nonumber\\
&\hspace{.2cm}=V_k^+(t)  \la(t,k)^* S(t,k)^{-1}(\a - {i} I_n)^{-1}
+V_k^-(t)\la(t,k)^* S(t,k)^{-1}(\a +{i} I_n)^{-1}.
\end{align*}
Here
\begin{align} 
&V_k^+(t) = I_{m}+j-{i} \la(t,k)^*S(t,k)^{-1}(\a - {i} I_n)^{-1}\la(t,k)j +\nn \\ 
&\hspace{.6cm} + {i} j \la(t,k)^*(\a^* - {i} I_n)^{-1}  S(t,k)^{-1}\la(t,k) + \la(t,k)^*S(t,k)^{-1} \times \nn \\ 
&\hspace{.6cm} \times (\a - {i} I_n)^{-1}\la(t,k)j \la(t,k)^*(\a^*
- {i} I_n)^{-1}S(t,k)^{-1}\la(t,k), \label{m34} 
\end{align}
and
\begin{align}
&V_k^-(t)=I_{m}-j+{i} \la(t,k)^*S(t,k)^{-1}(\a + {i} I_n)^{-1}\la(t,k)j+\nn \\
& - {i} j \la(t,k)^*(\a^* + {i} I_n)^{-1}S(t,k)^{-1}\la(t,k) -
\la(t,k)^*S(t,k)^{-1}\times \nn \\ 
& \times (\a + {i} I_n)^{-1}\la(t,k)j \la(t,k)^*(\a^* + {i} I_n)^{-1}S(t,k)^{-1}\la(t,k).  \label{m35}
\end{align}
It remains to show that $V_k^+(t)=H_k^+(t)$ and  $V_k^-(t)=H_k^-(t)$. To do this note that  \eqref{defHp2} and \eqref{defHm2}  tell us that    $H_k^+(t)$ and  $H_k^-(t)$ are also given by 
\begin{align} 
& H_k^+(t)= I_m+W_{\si_{t,k}}(i)jW_{\si_{t,k}}(-i)^*,\label{m36} \\[.1cm]
&H_k^-(t)=I_m-W_{\si_{t,k}}(-i)jW_{\si_{t,k}}(i)^*. \label{m37}
\end{align}
Now use \eqref{defWtk} and compute  the products. This shows  that $H_k^+(t)$ is equal to the right-hand side of \eqref{m34} and $H_k^-(t)$ is equal to the right-hand side of \eqref{m35}. In other words, $V_k^+(t)=H_k^+(t)$ and  $V_k^-(t)=H_k^-(t)$, and hence  \eqref{m28} is proved.

\medskip
\noindent\textsc{The fourth derivative in \eqref{dervas}.}   In this part we shall show that
\begin{align}
\left(\frac{d}{dt}W_{\si_{t, k}}\right)(\lambda)&=\left(\frac{H_k^+(t)}{\lambda-{i}}+\frac{H_k^-(t)}{\lambda+{i}}\right)
W_{\si_{t, k}}(\lambda)+\nn \\
&\hspace{2cm}
-2W_{\si_{t, k}}(\lambda)\left(\frac{P_1}{\lambda-{i}}+\frac{P_2}{\lambda+{i}}\right).
 \label{m43}
\end{align}
First, recall that $W_{\si_{t, k}}$ is given by 
\begin{equation}\label{defWtk}
W_{\si_{t, k}}(z)=I_m+i\la(t,k)^*S(t,k)^{-1}(zI_n-\a)^{-1}\la(t,k).
\end{equation}
Using \eqref{m18} and \eqref{m28}, we obtain
\begin{align} 
&\left(\frac{d}{dt}W_{\si_{t, k}}\right)(\lambda)={i} \Big(H_k^+(t)\la(t,k)^*S(t,k)^{-1}(\a-{i} I_n)^{-1}(\lambda I_n-\a)^{-1}\la(t,k)+\nn\\
&\hspace{1cm} +H_k^-(t)\la(t,k)^*S(t,k)^{-1}(\a+{i} I_n)^{-1}(\lambda I_n-\a)^{-1}\la(t,k)+\nn\\
&\hspace{1cm}  -2\la(t,k)^*S(t,k)^{-1}(\a-{i} I_n)^{-1}(\lambda  I_n-\a)^{-1}\la(t,k)P_1 +\nn\\
&\hspace{1cm} 
-2\la(t,k)^*S(t,k)^{-1}(\a+{i} I_n)^{-1}(\lambda I_n-\a)^{-1}\la(t,k)P_2 \Big).  \label{m38}
\end{align}
Clearly we have
\begin{align} 
& (\a-{i} I_n)^{-1}(\lambda  I_n-\a)^{-1}=(\lambda  -{i} )^{-1}\big((\a-{i} I_n)^{-1}+(\lambda  I_n-\a)^{-1}\big), \label{m39} \\
&  (\a+{i} I_n)^{-1}(\lambda  I_n-\a)^{-1}=(\lambda  +{i} )^{-1}\big((\a+{i} I_n)^{-1}+(\lambda  I_n-\a)^{-1}\big).  \label{m40}
\end{align}
Taking into account \eqref{defWtk}, \eqref{m39} and \eqref{m40}, we rewrite \eqref{m38} as
\begin{align} 
&\left(\frac{d}{dt}W_{\si_{t, k}}\right)(\lambda)=(\lambda -{i})^{-1}H_k^+(t)
\big( W_{\si_{t, k}}(\lambda)-W_{\si_{t, k}}(i)\big)+\nn\\
&\hspace{2cm} +(\lambda +{i})^{-1}H_k^-(t)
\big(W_{\si_{t, k}}(\lambda)-W_{\si_{t, k}}(-i)\big)+\nn\\
&\hspace{2cm} -2(\lambda -{i})^{-1} \big(W_{\si_{t, k}}(\lambda) -W_{\si_{t, k}}(i)\big)P_1 \nn\\
&\hspace{2cm} -2(\lambda +{i})^{-1}
\big(W_{\si_{t, k}}(\lambda)-W_{\si_{t, k}}(-i)\big)P_2. \label{m41}
\end{align}
Applying \eqref{m42p} and \eqref{m42m} with $\si_{t, k}$ in place of $\si_k$ we see that
\begin{equation}  \label{m42}
H_k^+(t)W_{\si_{t, k}}(i)=2W_{\si_{t, k}}(i)P_1, \quad H_k^-(t)W_{\si_{t, k}}(-i)=2W_{\si_{t, k}}(-i)P_2 .
\end{equation}
In view of \eqref{m42}, formula \eqref{m41} can be rewritten as \eqref{m43}.

\medskip
\noindent\textbf{Proof of Proposition \ref{propComp2}.} Let $Y$ be the matrix fucntion defined by \eqref{m44}.  The fact that $Y_k(t,z)$ is a solution of   equation \eqref{m5} follows by applying Theorem \ref{TmD1} with $\si_{t}$ in place of $\si$.   Similarly, using  \eqref{m43} and the fact that the matrices $P_1$, $P_2$, and $j$ are mutually commutative, we see that $Y_k(t,z)$ is a solution of  the system  \eqref{m6}.  Indeed,
\begin{align*}
\frac{d}{dt} Y_k(t,z)& =\left( \frac{d}{dt}W_{\si_{t, k}}\right)(-z) \left(I_m+\frac{{i}}{z}j\right)^k\exp\left\{-2t\left(\frac{P_1}{z+{i}}+\frac{P_2}{z-{i}}\right)\right\}  +  \\
&\hspace{3cm}-2Y_k(t,z)\left(\frac{P_1}{z+{i}}+\frac{P_2}{z-{i}}\right)\\
&=\left(\frac{H_k^+(t)}{-z-{i}}+\frac{H_k^-(t)}{-z+{i}}\right)Y_k(t,z)-2Y_k(t,z)\left(\frac{P_1}{-z-{i}}+\frac{P_2}{-z+{i}}\right)+\\
&\hspace{3cm}-2Y_k(t,z)\left(\frac{P_1}{z+{i}}+\frac{P_2}{z-{i}}\right)\\
&=-\left(\frac{H_k^+(t)}{z+{i}}+\frac{H_k^-(t)}{z-{i}}\right)Y_k(t,z).
\end{align*}
Thus $Y$ has the desired properties. \epr

\medskip
\begin{Rk}\label{RkDGHM2}
The equality $H_k^-=2I_m-(H_k^+)^*$ is immediate from \eqref{m36} and \eqref{m37}.
\end{Rk}


\medskip
\noindent\textbf{A remark on the continuous {space variable $x$} analogue.}  
{The transfer matrix function of the form \eqref{7.5} depending on the continuous parameter $t$ and discrete parameter $k$, which we use in this section, appeared in \cite{KaS},
whereas the transfer matrix function (of the form \eqref{7.5}) depending on continuous parameters $x$ and $t$ appeared first in \cite{ALS93, ALS94} (see also \cite{GKS1} and \cite{GKS2}).
To explain the case of    two continuous parameters in more detail  we begin with the following analogue of Lemma~\ref{lemdiscr1a}.  }
 
\begin{lem}\label{lemquadr2} Let  $\si=\{ \a, S_0,  \vt_1, \vt_2\}$   be a strongly   admissible quadruple, and let $k$ be a positive integer. Define
\begin{align}
&\vt_1(t)= e^{-it\a^k}\vt_1, \quad  \vt_2(t)= e^{it\a^k}\vt_2, \quad \la(t)=\begin{bmatrix} \vt_1(t) & \vt_2(t)\end{bmatrix},\nonumber\\
& S(t)= S_0+\int_0^t \sum_{\nu=1}^k  \a^{\nu-1}\la (s)j \la (s) (\a^*)^{k-\nu}\, ds, \quad t\geq 0.\label{defS(t)}
\end{align}
Then for each $t\in \BR$ the quadruple $\si(t)=\{ \a, S(t),  \vt_1(t), \vt_2(t)\}$   is  strongly admissible. 
\end{lem} 

 The proof of  the above lemma follows the same line  of reasoning as used in the proof of   Lemma \ref{lemdiscr1a}. We omit the details. 

Next, using the above lemma, we  construct a family of  pseudo-exponential potentials depending continuously on an additional   variable $t\in \BR$. The starting point is the   pseudo-exponential potential $v$ in \eqref{7.1} defined by the strongly admissible quadruple  $\si=\{ \a, S_0,  \vt_1, \vt_2\}$.   Fix $t\in \BR$, and let  $\si(t)$ be the strongly admissible quadruple  defined in Lemma \ref{lemquadr2}, i.e.,  $\si(t)=\{ \a, S(t),  \vt_1(t), \vt_2(t)\}$.  In particular, $S(t)$ is given by \eqref{defS(t)}.  We apply Lemma~\ref{lemadm1} with $\si(t)$ in place of  $\si$. Put
\begin{align}
\vt_{1}(x,t)&=e^{-ix\a}\vt_1(t)=e^{-i(x\a+it\a^k)}\vt_1,\nonumber \\
\vt_{2}(x,t)&=e^{ix\a}\vt_2(t)=e^{i(x\a+it\a^k)}\vt_2,\nonumber \\
\la(x,t)&= \begin{bmatrix} \vt_{1}(x,t) &\vt_{2}(x,t)\end{bmatrix}, \nonumber\\
 S(x,t)&= S(t)  +\int_0^x   \la (s, t)j \la (s, t)^* \, ds, \quad x\geq 0. \label{defSN(x,t)}
\end{align}
Then, by Lemma \ref{lemadm1},  the quadruple 
\begin{equation}
\label{defSiN} 
\si(x, t)=\{ \a, S(x,t),  \vt_{1}(x, t), \vt_{2}(x, t)\}
\end{equation}
is admissible for each $x\geq 0$. Furthermore,   the corresponding pseudo-exponential is given by
\begin{equation}\label{fampoten1}
v(x,t)= 2\vt_1^*e^{i(x\a^*+it(\a^*)^k)}S(x, t)^{-1}e^{i(x\a+it\a^k)}\vt_{2}, \quad x\geq 0.
\end{equation}
In this way we obtain a family of  pseudo-exponential potentials depending on the additional parameter {$t$}. 

For {$k=2,3$} the potential $v$ in \eqref{fampoten1} satisfies {important integrable} nonlinear equations; see, e.g.,  \cite[Section 2]{ALS93} and   \cite[Section 4]{GKS2}.  {Namely}, for $k=2$ the function $v$ is a  solution of the matrix nonlinear Schr\"odinger equation  
\begin{equation}\label{GF2}
2 \frac{ \partial v}{ \partial t}+i\frac{ \partial ^{2} v}{ \partial x^{2}}
+2ivv^{*}v=0,
\end{equation}
 and for $k=3$ it satisfies the matrix modified Korteweg-de Vries equation 
 \begin{equation}
4\frac{ \partial v}{ \partial t}+\frac{ \partial ^{3}v}{ \partial x^{3}}
+3\Big(\frac{ \partial v}{ \partial x}v^{*}v+vv^{*}\frac{ \partial v}{ \partial x}\Big)
=0.
\end{equation}

The proofs of the above results can be obtained by direct computations as in \cite[Section 2]{ALS93} and   \cite[Section 4]{GKS2} or by  {applying} the generalized  B\"acklund-Darboux transformation \cite{SaAJMAA}  {to auxiliary linear systems and using compatibility condition (zero curvature equation)
$G_t-F_x+GF-FG=0$, which is a continuous case equivalent of  the compatibility condition} \eqref{m4}.
For instance, equation \eqref{GF2} is equivalent to the compatibility condition  $G_t-F_x+GF-FG=0$, where
$G=i z j +V(x,t)$, 
$$F=i\big(z^2 j-izjV(x,t)-\big(V_x(x,t)+jV(x,t)^2)/2\big),$$ 
and $V$ has the form \eqref{se2}.

\appendix
\section{Appendix: admissible quadruples}\label{App}
\renewcommand{\theequation}{A.\arabic{equation}}
\setcounter{equation}{0}
Let {$S_0$} and $\a$ be $n\times n$ matrices, and let $\vt_1$ and $\vt_2$ be matrices of sizes $n\ts m_1$ and $n\ts m_2$, respectively. The quadruple $\{ \a, S_0,  \vt_1, \vt_2\}$   {is called  \emph{admissible} (see the beginning of Section \ref{dies})} if 
\begin{equation}
\label{defquadr}
S_0 > 0 \ands  \a S_0-S_0 \a^*=i \la \la^*,  \ \mbox{where  $\la:=
\begin{bmatrix}\vt_1&\vt_2 \end{bmatrix}$}.
\end{equation}
{By definition the \emph{transfer  function}} {associated with} the quadruple 
$$\si=\{ \a, S_0,  \vt_1, \vt_2\}$$ 
is the matrix function $W_{\si}$ given by  (see \eqref{7.5} and \eqref{d7}):
\begin{align}
W_\si(z)&= I_m +i \la^*S_0^{-1}(z I_n- \a)^{-1}\la, \  \mbox{where $m:=m_1+m_2$}.\label{defW}
\end{align} 
{With $\si$ we associate  two other} rational matrix functions:
\begin{align}
{\va_{1,\si}}(z)&=i\vt_2^*S_0^{-1}(zI_n-\beta_1)^{-1}\vt_1, \  \mbox{where $\b_1:=\a -i\vt_1\vt_1^* S_0^{-1}$};\label{defWeyl1} \\
{\va_{2,\si}}(z)&=-{i} \vartheta_{1}^{*}S_0^{-1}(z  I_{n} + \b_2)^{-1}
\vartheta_{2}, \quad \mbox{where $\beta_2:=\a-{i} \vt_2\vt_2^*S_0^{-1}$}.
\label{defWeyl2}
\end{align} 
{If   $\si=\{ \a, S_0,  \vt_1, \vt_2\}$ is an admissible quadruple, then the same holds true for the} 
{associate} {quadruple 
\begin{align}
\si^\#:=\{ \a, S_0,  \vt_2, \vt_1\}
\label{va12!}
\end{align} 
and a direct computation shows that }
\begin{equation}\label{va12}
\va_{1, \si^\#}(-z)=\va_{2, \si}(z). 
\end{equation}

{Admissible quadruples are closely related to symmetric S-nodes (see \cite[Section 2.2]{SaL3}). Indeed, if}   $\si=\{ \a, S_0,  \vt_1, \vt_2\}$ is an admissible quadruple, then  the  matrices 
\[
 {A:=\a, \quad S:=S_0, \quad  \Pi_1:=\la, \quad  \ga_1:=S_0^{-1}\la}, 
\]
form a  symmetric S-node  {(with $J=I_n$) as defined in  {Section 2.2} }of  \cite{SaL3}. Using this connection, it is readily seen that  item (iii)  in the following lemma is a special case of  Theorem 2.2.1 in \cite{SaL3} (see also  {\cite{SaL1}}, Theorem 17.1 in \cite{BGKR2}  {or Corollary 1.15 in \cite{SaSaR}}).  

\begin{lem}\label{lem11} Let  $\si=\{ \a, S_0,  \vt_1, \vt_2\}$ be an admissible quadruple. Then
\begin{itemize}
\item[\textup{(i)}]  $\s(\a)\subset \overline{ \BC_+}$;
\item[\textup{(ii)}] {if, in addition, the pair $\{\a, \la\}$ is controllable, then  $\s(\a)\subset  { \BC_+}$, and in that case $S_0$ is uniquely determined by $\a$ and  $\la$.
 {More precisely,} $S_0$  is given by
\begin{equation}
\label{altS0}
S_0=\frac{1}{2\pi}\int_{-\iy}^{\iy} (\lambda  I_n- \a)^{-1}\la \la ^*(\lambda  I_n-\a^*)^{-1}d\lambda. 
\end{equation}}
\item[\textup{(iii)}]  $W_\si$ is a rational Schur class function and its values on the real line are unitary matrices.  {Moreover, we have }
\begin{align}
& { I_m-W_{\si}(z)W_{\si}(\ov{\zeta} )^*=}\nn \\
&\hspace{1cm} {=i (z-\zeta )\la^*S_0^{-1}(zI_n-\a)^{-1}S_0(\zeta I_n-\a)^{-1}S_0^{-1}\la,} \label{basicineq1}\\
& I_m-W_\si(z)^*W_\si(z)=i(z-\bar{z})\la^*(\bar{z}I_n-\a^*)^{-1}S_0^{-1}(zI_n-\a )^{-1}\la. \label{basicineq2} 
\end{align}
\end{itemize}
\end{lem}
\bpr In remains to prove items (i) and (ii). Since $S_0$ is positive definite, the second part of \eqref{defquadr} implies that 
\[
 { \left(S_0^{-1/2} \a  S_0^{1/2}\right)- \left(S_0^{-1/2} \a  S_0^{-1/2}\right)^*=iS_0^{-1/2}\la \la^*S_0^{-1/2},  \quad  S_0^{-1/2}\la \la^*S_0^{-1/2}\geq 0.}
\]
It follows that the numerical range  of  the matrix $S_0^{-1/2} \a  S_0^{1/2}$ is a subset of  $\overline{ \BC_+}$, and hence $\s(\a)=\s(S_0^{-1/2} \a  S_0^{1/2})\subset \overline{ \BC_+}$.   This proves item (i).  

 {The first statement in} item (ii) is a  straightforward application of the classical Chen-Wimmer inertia theorem (see Corollary 1 in \cite[Section 13.1]{LT85}).
{For  \eqref{altS0} see  \cite{IvSa} or Theorem I.4.1 in \cite{GGK1} (see, also  \eqref{d23}).}

\begin{lem}\label{lem12} Let  $\si=\{ \a, S_0,  \vt_1, \vt_2\}$ be an admissible quadruple, and partition  {its transfer function} $W_\si$ as a $2\ts 2$ block matrix function,
\begin{equation}\label{d26}
W_\si(z)= \begin{bmatrix} a(z) & b(z ) \\[.2cm] c(z ) & d(z ) \end{bmatrix}, \ \mbox{where $a(z)$ has size  $m_1 \times m_1$}.
\end{equation}
Then the {associated  functions  $\va_{1,\si}$ and $\va_{2,\si}$ are given by  
\begin{equation}\label{abcd}
\va_{1,\si}(z)=c(z)a(z)^{-1} \ands  \va_{2,\si}(z)=b(-z)d(-z)^{-1}.
\end{equation}}
\end{lem}
\bpr  From \eqref{defW}  and $\la=\begin{bmatrix}  \vt_1 & \vt_2 \end{bmatrix}$ it follows that
\[
a(z)  = I_{m_1}+i\vt_1^*S_0^{-1}(z I_n- \a)^{-1}\vt_1,  \quad  c(z)   = i\vt_2^*S_0^{-1}(z I_n- \a)^{-1}\vt_1.
\]
Since  $a$ is a proper rational matrix function whose value  at infinity is the $m_1\ts m_1$ identity matrix, a classical inversion theorem (see, e.g., Theorem 2.2 in   \cite{BGKR1})  tells us that
\[
a(z)^{-1}= I_{m_1}-i\vt_1^*S_0^{-1}(z I_n- \a^\ts )^{-1}\vt_1,  \ \mbox{{where $\a^\ts=\a -i\vt_1\vt_1^* S_0^{-1}=\b_1$}}. 
 \]
 { Hence, again using Theorem 2.2 in   \cite{BGKR1} and taking into account \eqref{defWeyl1}, we know that 
\begin{align*}
c(z) a(z)^{-1}&= i\vt_2^*S_0^{-1}(z I_n- \a)^{-1}\vt_1a(z)^{-1}
=i\vt_2^*S_0^{-1}(z I_n- \a^\ts)^{-1}\vt_1\\
&{= i\vt_2^*S_0^{-1}(z I_n- \b_1)^{-1}\vt_1=\va_{1, \si}(z),}
\end{align*}}
which proves the first part of \eqref{abcd}.

{To prove the second  identity in   \eqref{abcd} we use the associate admissible quadruple $\si^\#=\{ \a, S_0,  \vt_2, \vt_1\}$ and  equality \eqref{va12}. Note that 
\begin{equation}\label{Ad1}
\begin{bmatrix} {\vt}_2&{\vt}_1 \end{bmatrix}=\begin{bmatrix} \vt_1& \vt_2 \end{bmatrix}U, \   \mbox{where   $U=\begin{bmatrix} 0 & I_{m_1} \\ I_{m_2}& 0 \end{bmatrix}$}.
\end{equation}
Since $U$ is unitary,   it is clear from \eqref{Ad1} and \eqref{defW} that $W_{\si^\#}=U^* W_\si U$. But then, using the special form of $U$, we obtain 
\begin{equation}\label{Ad2}
W_{\si^\#}=U^* W_\si U=  \begin{bmatrix} d(z) & c(z ) \\[.2cm] b(z ) & a(z ) \end{bmatrix}.
\end{equation}
Applying the result of the first part of the proof  to $\si^\#$ in place of  $\si$ we see  that $\va_{1, \si^\#}=b(z)d(z)^{-1}$. But then   \eqref{va12} yields the second identity in \eqref{abcd}}.\epr

\begin{thm}\label{thminv} Let $\va$ be an  $m_2 \ts m_1$ rational matrix function. Then  $\va$   {coincides {with} some   function $\va_{1,\si}$ associated with {an} admissible quadruple 
$\si$ if and   only if  $\va$ is strictly proper. In that case, the corresponding $\si$} can be obtained explicitly by using the  following procedure.
\begin{description}
\item[\textup{Step 1.}] Let $n$  be the McMillan degree of $\va$, and construct a minimal realization of $\va:$
\begin{equation}\label{c15} 
\va(z)=i \th_2^*(z I_n
- \g)^{-1}
\th_1.
\end{equation}
\item[\textup{Step 2.}]
Choose {$X$ to be the unique} positive definite solution   of the algebraic Riccati equation
\begin{equation}
\label{Ricceq1} 
\g X-X\g^*-iX\th_2\th_2^*X +i \th_1\th_1^*=0.
\end{equation}
\item[\textup{Step 3.}]  Put 
\begin{align} 
&S_0=I_n, \quad \vt_{1}=X^{ -{1}/{2}} \theta_{1}, \quad \vt_{2}=X^{{1}/{2}}  \theta_{2}, \label{c17a}\\[.2cm]
&\hspace{1.4cm} \a=X^{- {1}/{2}}\g  X^{ {1}/{2}}+i \vt_1\vt_1^*. \label{c17b}
\end{align} 
\end{description}
With these  choices the quadruple $\si= \{\a, S_0, \vt_1,  \vt_2\}$ is   admissible   and $\va$  {coincides with its {first  associate function $\va_{1, \si}$}}. Moreover, in this case the pairs $\{\a, \vt_1\}$ and 
$\{\a, \vt_2\}$ are controllable.
\end{thm}
\bpr  First note that the minimality of the realization \eqref{c15} is equivalent to the requirement that simultaneously the pair $\{\th_2^*, \g\}$ is observable and the pair $\{\g, \th_1\}$ is controllable.  But then we can apply Proposition 2.2 in \cite{GKS2} to show that the Riccati equation \eqref{Ricceq1} has a  positive definite solution $X$, a result which has its roots in Kalman's theory of  mathematical systems \cite{KFA} (see also 
 \cite[pp. 358 and 369]{LR},  {where the uniqueness of $X$ is shown as well}). The remaining part of the proof is split into two parts.

\smallskip \noindent
\textsc{Part 1.} Using the definition of $\a$ in \eqref{c17b}, the {unique  positive definite solution $X$} of the Riccati equation \eqref{Ricceq1}, and the definitions of  $\vt_1$ and $\vt_2$ in \eqref{c17a} we see that
\begin{align*} 
\a-\a^* &=X^{-1/2}\g  X^{1/2}-X^{1/2}\g^*X^{- 1/2}+2i \vt_1\vt_1^*\\[.1cm]
&=iX^{1/2}\th_2\th_2^*X^{1/2}-iX^{-1/2}\th_1\th_1^*X^{-1/2}+2i \vt_1\vt_1^*\\[.1cm]
&=i\vt_2\vt_2^*- i \vt_1\vt_1^* +2i \vt_1\vt_1^*
=i\vt_2\vt_2^*+ i \vt_1\vt_1^*.
\end{align*} 
Since $S_0=I_n$ by definition, we conclude that $\si= \{\a, S_0, \vt_1,  \vt_2\}$ is   an admissible quadruple. Moreover, the associate  {function  $\va_{1, \si}$} is equal to $\va$. Indeed, in this case using \eqref{c17b}, we have  
\begin{align} \label{Ad7}
{\b_1}=\a- i\vt_1\vt_1^*= X^{-1/2}\g X^{1/2}.
\end{align}
Thus { $X^{1/2}\b_1 X^{-1/2}=\g$}, and hence 
{\begin{align*} 
\va_{1,\si}(z)&= i\vt_2^*S_0^{-1}(zI_n-\b_1)^{-1}\vt_1\\[.1cm]
&=i\th_2^*X^{1/2}(zI_n-\b_1)^{-1}X^{-1/2}\th_1\\[.1cm]
&=i\th_2^*(zI_n-X^{1/2}\b_1 X^{-1/2})^{-1}\th_1\\[.1cm]
&=i\th_2^*(zI_n-\g)^{-1}\th_1=\va.
\end{align*} }

\smallskip \noindent
\textsc{Part 2.} In this part we prove that the pairs $\{\a, \vt_1\}$ and  $\{\a, \vt_2\}$ are controllable. 
 {Since the realization \eqref{c15} is minimal and $\va={\va_{1,\si}}$, the same is valid for the realization \eqref{defWeyl1}, and hence the pairs {$\{\b_1, \vt_1\}$ and $\{ (\b_1)^*, \vt_2 \}$ }} are controllable.   {Formula \eqref{defWeyl1}   implies  that  { $\a=\b_1+i\vt_1\vt_1^*$}. Using the latter identity, the fact that {$\{ \b_1, \vt_1\}$}} is controllable and the equality $\im \vt_1= \im \vt_1\vt_1^*$, we conclude that $\{\a, \vt_1\}$ is controllable too. 

 It remains to prove that $\{\a, \vt_2\}$ is controllable. Since $S_0=I_n$ and   the quadruple $\si= \{\a, S_0, \vt_1,  \vt_2\}$ is  admissible, we have $\a - \a^*=i(\vt_1\vt_1^*+\vt_2\vt_2^*)$, and therefore 
 \[
{ \b_1}=\a-i\vt_1\vt_1^*=\a^*+ i\vt_2\vt_2^*=(\a - i\vt_2\vt_2^*)^*.
 \]
Thus $\a= {(\b_1)}^*+i\vt_2\vt_2^*$. We know that  $\{{(\b_1)}^*, \vt_2\}$ is controllable and that $\im \vt_2= \im \vt_2\vt_2^*$.  {These  facts} imply that  $\{\a, \vt_2\}$ is  also controllable.\epr

\medskip
Using the identity \eqref{va12} and Theorem \ref{thminv} we can show that a proper rational matrix function $\va$ also coincides with the   second associate  function $\va_{2,\si}$ of some admissible quadruple $\si$. More precisely, we have the following corollary.

\begin{Cy} \label{CyA5}  Let $\wt{\va}$ be an  $m_2 \ts m_1$   rational matrix function. Then  $\wt{\va}$  coincides  with  some   function $\va_{2,\si}$ associated with {an} admissible quadruple 
$\si$ if and   only if  $\va$ is strictly proper.  The corresponding $ \si= \{\a, S_0,  \vt_1,   \vt_2\}$ can be obtained explicitly by using the  following procedure.
First, we construct a minimal realization of $\wt \va$,
\[
\wt{\va}(z)=-i\th_1^*(zI_n+\g)^{-1}\th_2.
\] 
Next, we choose $X$ to be the unique positive definite solution of the Riccati equation
$\g X-X\g^*-iX\th_1\th_1^*X +i \th_2\th_2^*=0$. Finally, we put
\begin{align}
& S_0=I_n, \quad \vt_{1}=X^{ {1}/{2}} \theta_{1}, \quad \vt_{2}=X^{-{1}/{2}}  \theta_{2},   \label{Ad91}\\
& \hspace{1cm} \a=X^{- {1}/{2}}\g  X^{ {1}/{2}}+i \vt_2\vt_2^*.  \label{Ad92}
\end{align} 
Then $\wt{\va}= \va_{2,\si}$  and the pairs $\{\a, \vt_1\}$ and  $\{\a, \vt_2\}$ are controllable.
\end{Cy}
\bpr Put $\va(z)=\wt{\va}(-z)$. Then $\va$ is rational, and  $\wt{\va}$ is  strictly proper if and only if   $\va$ is strictly proper. But,  by   Theorem \ref{thminv},     the latter happens if and only if   there exists an admissible quadruple $\wt{\si}= \{\a, S_0, \wt \vt_1,  \wt \vt_2\}$ such that $ \va(z)=\va_{1, \wt{\si}}(z)$. Next, {taking into account \eqref{va12!}, we see that
$\wt{\si}$ is the associate quadruple $\si^\#$ for
${\si}=  \{\a, S_0, \wt{\vt}_2,  \wt{\vt}_1\} $}. Then we know from  \eqref{va12} that
\[
{\wt{\va}(z)=\va(-z)=\va_{1,\si^{\#}}(-z)=\va_{2,\si}(z).}
\]
 Moreover, Theorem \ref{thminv} provides a method to construct $\wt{\si}$, which in terms of  {$\wt \si =\si^{\#}$ (or, equivalently, $\si= \wt{\si}^\#$) }
 yields the procedure  to recover $ \si$ from $\wt{\va}$  and 
 formula \eqref{defWeyl2}.\epr

\medskip\noindent

 {{\bf Proof of Proposition \ref{propdefv}}.}
{  Let $S(x)$ be the $n\ts n$ matrix defined  {by  \eqref{7.2}}, and put
\[
 Q(x)=\begin{bmatrix}S_0 & -iI_n\end{bmatrix}e^{-2ixA}\begin{bmatrix}I_n \\ 0\end{bmatrix}, \quad  x\in \BR_+,
\]
 {where $A$ is given by \eqref{defA}. In order to prove Proposition \ref{propdefv}}, that is, to prove \eqref{altdefv},} 
it suffices to show that   $Q(x)= \de(x)$, where $\de(x)=e^{-ix\a} S(x)e^{-ix\a^*}$, for each $x\in \BR_+$.   {For that  we  directly
differentiate $\de(x)$ using the first identity  in \eqref{defSx} (where we substitute
the expression for $\la$ from \eqref{7.2})  and obtain
\begin{align}
&\frac{d}{dx}\de(x)=-i\big(\a \de(x)+\de(x)\a^*\big) + e^{-2ix\a} \vt_1\vt_1^*-\vt_2\vt_2^* e^{-2ix\a^*}\label{derDe}.
\end{align}
Now, we show that $Q$ satisfies the same first order linear differential equation }
\begin{align}
&\frac{d}{dx}Q(x)=-i\big(a Q(x)+Q(x)\a^*\big) + e^{-2ix\a} \vt_1\vt_1^*-\vt_2\vt_2^* e^{-2ix\a^*}.\label{derQ}
\end{align}
To prove \eqref{derQ} note that 
\[
i\frac{d}{dx}Q(x)= \begin{bmatrix} S_0 & -iI_n \end{bmatrix}\Big( A e^{-2ixA}+e^{-2ixA}A\Big)  \begin{bmatrix} I_n \\ 0\end{bmatrix}.
\]
Using \eqref{7.3} and the definition of $A$ in \eqref{defA} we see that 
\begin{align*}
\begin{bmatrix} S_0 & -iI_n \end{bmatrix}A&=\begin{bmatrix}S_0\a^*  + i\vt_1\vt_1^* & -i\a \end{bmatrix}\\
&=\begin{bmatrix}\a S_0   - i\vt_2\vt_2^* & -i\a \end{bmatrix}\\
&=\a \begin{bmatrix} S_0 & -iI_n \end{bmatrix} -i \begin{bmatrix} \vt_2\vt_2^*  & 0 \end{bmatrix}.
\end{align*}
Furthermore, we have 
\[
A\begin{bmatrix} I_n \\ 0\end{bmatrix}= \begin{bmatrix} \a^*\\ -\vt_1\vt_1^* \end{bmatrix}= \begin{bmatrix}I_n\\ 0 \end{bmatrix} \a^* + \begin{bmatrix}0\\ I_n \end{bmatrix}(-\vt_1\vt_1^*).
\]
It follows (using again the definiton of $A$ in \eqref{defA})  that
\begin{align}\nn
i\frac{d}{dx}Q(x)&=\a \begin{bmatrix} S_0 & -iI_n \end{bmatrix}e^{-2ixA} \begin{bmatrix}I_n\\ 0 \end{bmatrix}  -i \begin{bmatrix} \vt_2\vt_2^*  & 0\end{bmatrix}e^{-2ixA}  \begin{bmatrix}I_n\\ 0 \end{bmatrix}  \\ \nn
&\hspace{2cm}+\begin{bmatrix} S_0 & -iI_n \end{bmatrix}e^{-2ixA} \begin{bmatrix}I_n\\0\end{bmatrix}\a^*\\ \nn
&\hspace{2cm}+ \begin{bmatrix} S_0 & -iI_n \end{bmatrix}e^{-2ixA}  \begin{bmatrix}0\\ I_n \end{bmatrix}(-\vt_1\vt_1^* )\\ \label{eqQ}
&=\a Q(x)-i\vt_2\vt_2^* e^{-2ix\a^*}+Q(x)\a^*+i e^{-2ix\a} \vt_1\vt_1^*.
\end{align}
 {Equation \eqref{derQ} is immediate from $Q$. Since $\de$ and $Q$ satisfy the same first order linear differential equation with the same
initial condition $\de(0)=Q(0)=S_0$, we see that} $\de(x)=Q(x)$ for each $x\in \BR_+$, as desired.

\bigskip 
\noindent{\bf Acknowledgments.}
 {The research of A.L. Sakhnovich was supported by the
Austrian Science Fund (FWF) under Grant  No. P24301.}
The authors are grateful to I.~Roitberg for her help and useful remarks.


\end{document}